\documentclass[11pt,a4paper]{article}
\usepackage[T1]{fontenc}
\usepackage[utf8]{inputenc}
\usepackage{lmodern}
\usepackage{amsmath,amssymb,amsthm,mathtools}
\usepackage[margin=27mm,headheight=14pt]{geometry}
\usepackage{microtype}
\usepackage{enumitem}
\usepackage{xcolor}
\usepackage{fancyhdr}
\usepackage{hyperref}
\hypersetup{colorlinks=true,linkcolor=blue,citecolor=blue,urlcolor=blue,
	pdftitle={Blocks with one simple module and solvable groups},
	pdfauthor={Kun Zhang and Yuanyang Zhou},
	pdfsubject={Modular representation theory of finite groups},
	pdfkeywords={block algebra, Morita equivalence, nilpotent block}}
\newtheorem{theorem}{Theorem}[section]
\newtheorem{proposition}[theorem]{Proposition}
\newtheorem{lemma}[theorem]{Lemma}
\newtheorem{corollary}[theorem]{Corollary}
\theoremstyle{definition}
\newtheorem{definition}[theorem]{Definition}
\theoremstyle{remark}

\numberwithin{equation}{section}
\newcommand{\cO}{\mathcal{O}}
\newcommand{\Fp}{\mathbb{F}_p}
\newcommand{\Mor}{\mathrm{Mor}}
\DeclareMathOperator{\Irr}{Irr}
\DeclareMathOperator{\IBr}{IBr}
\DeclareMathOperator{\Aut}{Aut}
\DeclareMathOperator{\Inn}{Inn}
\DeclareMathOperator{\Out}{Out}
\DeclareMathOperator{\End}{End}
\DeclareMathOperator{\Hom}{Hom}
\DeclareMathOperator{\Br}{Br}
\DeclareMathOperator{\GL}{GL}

\DeclareMathOperator{\Syl}{Syl}

\DeclareMathOperator{\id}{id}
\DeclareMathOperator{\rank}{rank}

\newcommand{\op}{\mathrm{op}}
\newcommand{\proofstep}[1]{%
	\par\addvspace{\medskipamount}%
	\noindent\textbf{#1}\par\nobreak
	\smallskip\nobreak\noindent\ignorespaces
}
\setlist[enumerate]{label=\textup{(\roman*)},leftmargin=2em,itemsep=3pt,topsep=5pt}
\title{Blocks with one simple module\\and solvable groups}
\author{%
	Kun Zhang\\
	{\footnotesize Faculty of Mathematics and Statistics, Hubei University, Wuhan 430062, China.}\\[0.55em]
	Yuanyang Zhou\\
	{\footnotesize School of Mathematics and Statistics, Central China Normal University, Wuhan 430079, China.}%
}
\date{}
\begin{document}
	\maketitle
	\vspace{-2em}
	\begin{abstract}
		Let $b$ be a block of a finite group with one isomorphism class of simple modules in characteristic $p$.
		We construct, over a sufficiently large complete discrete valuation ring of characteristic zero, a Morita equivalence between $b$ and a block of a finite solvable group.
		The proof uses equivariant Morita bimodules for the quasisimple blocks that occur in the component reduction.
		We also prove that a block over an algebraically closed field of characteristic $p$ is nilpotent whenever it is Morita equivalent to the group algebra of a finite $p$-group.
		For the latter result, we establish a criterion for twisted group algebras using derivations of the basic algebra and group transfer.
	\end{abstract}

	\section{Introduction}\label{sec:introduction}
	
	For a block $b$, let $l(b)=|\IBr(b)|$ be the number of isomorphism classes of simple modules in characteristic $p$.
	Puig's structure theorem shows that a nilpotent block is Morita equivalent to the group algebra of a defect group and therefore has one simple module \cite[Theorem~1.6]{Puig1988}.
	The converse fails, even for solvable groups \cite[Theorem~3]{Sambale2024}.
	The condition $l(b)=1$ thus defines a broader class whose description up to Morita equivalence remains a natural structural problem.
	
	Kessar and Linckelmann conjectured that every block with one simple module is Morita equivalent over the coefficient ring to a block of a finite $p$-solvable group; see \cite[Introduction]{MalleNavarroSpath2018}.
	We prove that the group can be chosen solvable.
	We fix a sufficiently large splitting $p$-modular system $(K,\cO,k)$, where $\cO$ is a complete discrete valuation ring, $\operatorname{char}(K)=0$, and $k$ is algebraically closed of characteristic $p$.
	The precise conventions are given at the beginning of Section~\ref{sec:preliminaries}.
	
	\begin{theorem}\label{thm:solvable}
		Let $G$ be a finite group and let $b$ be a block of $\cO G$ with $l(b)=1$.
		There are a finite solvable group $H$ and a block $f$ of $\cO H$ such that $\cO Gb$ and $\cO Hf$ are Morita equivalent over $\cO$.
	\end{theorem}
	
	Our second result concerns recognition of nilpotent blocks from their modular Morita equivalence class.
	Puig proved the corresponding assertion for block algebras over $\cO$ \cite[Theorem~8.2]{Puig1999}.
	The following theorem only assumes a Morita equivalence over the residue field.
	
	\begin{theorem}\label{thm:nilpotency}
		Let $k$ be algebraically closed of characteristic $p$, let $G$ be a finite group, and let $b$ be a block of $kG$.
		If $kGb$ is Morita equivalent to $kP$ for a finite $p$-group $P$, then $b$ is nilpotent.
	\end{theorem}
	
	The proof of Theorem~\ref{thm:solvable} reduces to a precise equivariance statement for quasisimple groups.
	A finite group is quasisimple if it is perfect and its quotient by its center is nonabelian simple.
	We consider blocks for which the stabilizer $\Aut(X)_c$ acts transitively on $\IBr(c)$.
	For these quasisimple groups, the blocks in question have abelian defect groups and are inertial \cite[Propositions~1.2 and~4.1]{JiangZhangZhou2026}.
	For the reduction, we need to choose the Morita bimodules with compatible central and automorphism actions.
	
	\begin{proposition}\label{prop:local}
		Let $X$ be a quasisimple group, and let $c$ be a block of
		$\cO X$ with defect group $D$ such that $\Aut(X)_c$ acts
		transitively on $\IBr(c)$. Then $D$ is abelian, and $(X,c,D)$
		admits a scalar-equivariant Morita datum in the sense of
		Definition~\ref{def:datum}.
	\end{proposition}
	
	To prove Theorem~\ref{thm:nilpotency}, we apply Theorem~\ref{thm:solvable} to $\cO Gb$.
	Fong reduction over $\cO$ gives a twisted group algebra $\cO^\alpha L$ with $O_{p'}(L)=1$ and a cocycle taking values in a finite group of $p'$-roots of unity.
	We prove that Morita equivalence of $k^{\bar\alpha}L$ to a $p$-group algebra forces $L$ to be a $p$-group.
	The proof tests derivations on the first radical layer of the basic algebra. Group transfer and the action on the center of a Thompson critical subgroup and the count of projectively regular classes in characteristic two exclude a nontrivial quotient $L/O_p(L)$. The twisting then vanishes over $\cO$, and Puig's theorem for block algebras over $\cO$ applies.
	
	Section~\ref{sec:preliminaries} develops scalar-equivariant Morita data,
	the tensor construction, and the extension proposition
	(Proposition~\ref{prop:quotient-induction}).
	Section~\ref{sec:quasisimple} constructs the quasisimple data required in
	Proposition~\ref{prop:local}. Section~\ref{sec:reduction} combines these
	two inputs: the component proposition produces a new block, and a strict
	decrease in the index of the central $p'$-subgroup proves
	Theorem~\ref{thm:solvable}.
	Section~\ref{sec:derivations} proves the derivation criterion for twisted group algebras. Finally, Section~\ref{sec:integral} gives Fong reduction over
	$\cO$ and proves Theorem~\ref{thm:nilpotency}.

	\section{Equivariant Morita data and component extensions}\label{sec:preliminaries}
	
	This section records the constructions needed in the reduction. We first
	specify the equivariant data, then form tensor powers, and finally pass
	from an orbit of components to a block of a central extension of a
	normalizer. Section~\ref{sec:quasisimple} supplies the required
	quasisimple data.
	
	\subsection{Coefficients and block notation}
	
	Fix a prime $p$. Let $\cO$ be a complete discrete valuation ring with fraction
	field $K$ of characteristic zero and algebraically closed residue field $k$ of
	characteristic $p$. We assume that $K$ splits the finite groups occurring in the
	argument, including the groups introduced by the extension constructions. We
	keep $\cO$ fixed. Unless another ground ring is specified, tensor products are
	taken over $\cO$.
	
	All groups are finite and all modules are finitely generated. We identify a
	block with its primitive central idempotent and use the same letter for its
	reduction modulo the maximal ideal of $\cO$. For a block $b$ of $\cO G$, write
	$\Irr(b)$ for its irreducible ordinary characters and $\IBr(b)$ for its
	irreducible Brauer characters, and put $l(b)=|\IBr(b)|$. If $L\unlhd G$,
	then $b$ covers a block $c$ of $\cO L$ when $bc\ne0$. If $Z\leq Z(G)$ and
	$\pi:\cO G\longrightarrow\cO(G/Z)$ is the quotient map, a block $\bar b$
	of $\cO(G/Z)$ is dominated by $b$ when $\pi(b)\bar b=\bar b$.
	
	We write $O_p(G)$ and $O_{p'}(G)$ for the largest normal $p$-subgroup and
	$p'$-subgroup of $G$, respectively. For $x\in G$, let $c_x$ denote the inner
	automorphism $y\mapsto xyx^{-1}$. We use $c_x^H$ when the group $H$
	on which conjugation acts needs to be specified. Thus
	\[
	\Inn(G)=\{c_x:x\in G\},\qquad \Out(G)=\Aut(G)/\Inn(G).
	\]
	Subscripts denote stabilizers; for example,
	\[
	\Aut(G)_{D,b}=\{a\in\Aut(G):a(D)=D,\ a(b)=b\}.
	\]
	For a group $Q$ and an automorphism $a$ of $Q$, set
	$[Q,a]=\langle x^{-1}a(x):x\in Q\rangle$.
	
	For a $p$-subgroup $Q\leq G$, the Brauer map
	$\Br_Q:(\cO G)^Q\longrightarrow kC_G(Q)$ retains the coefficients of the
	elements centralizing $Q$ and reduces them modulo the maximal ideal. A
	$b$-Brauer pair is a pair $(Q,e)$ with $e$ a block of $kC_G(Q)$ such that
	$\Br_Q(b)e\ne0$. The maximal $b$-Brauer pairs form a single $G$-conjugacy
	class \cite[Theorem~3.10]{AlperinBroue1979}; their subgroup components are the defect groups
	of $b$. We also write $e$ for its lift to $\cO C_G(Q)$. Given a defect group
	$D$ of $b$, its Brauer correspondent in $N_G(D)$ is denoted by $b_0$.
	A block is nilpotent if $N_G(Q,e)/C_G(Q)$ is a $p$-group for every one of its
	Brauer pairs. We use basic Morita equivalence in the sense of
	Puig~\cite[Corollary~7.4]{Puig1999}: the inducing indecomposable bimodule has
	an endopermutation source of rank prime to $p$. A block is inertial if it is basically
	Morita equivalent to its Brauer correspondent \cite[Definition~2.16]{Puig2011}.
	
	\subsection{Central balance and semilinear maps}
	
	For an algebra $A$, let $Z(A)$, $J(A)$, and $A^{\op}$ denote its center,
	Jacobson radical, and opposite algebra. An $(A,A_0)$-bimodule $U$ is a Morita
	bimodule if $U\otimes_{A_0}-$ is an equivalence from left $A_0$-modules to left
	$A$-modules. An inverse Morita bimodule $U^\vee$ satisfies
	\begin{equation}\label{eq:morita-identities}
		U\otimes_{A_0} U^\vee\cong A,\qquad
		U^\vee\otimes_A U\cong A_0.
	\end{equation}
	We use the standard identifications
	\[
	\End_A(U)\cong A_0^{\op},\qquad \End_{A\text{-}A_0}(U)\cong Z(A)
	\]
	\cite[Theorem~2.8.2 and the proof of Corollary~2.8.6]{Linckelmann2018a}.
	The Morita bimodules used over $\cO$ are $\cO$-free. Reducing
	\eqref{eq:morita-identities} modulo the maximal ideal gives a Morita equivalence
	over $k$; in particular, an $\cO$-linear Morita equivalence between block
	algebras preserves the number of simple modules.
	
	An idempotent $e\in A$ is full if $AeA=A$, or equivalently if
	$1_A=\sum_{i=1}^m a_i e b_i$ for some $a_i,b_i\in A$.
	For full idempotents $e\in A$ and $f\in A_0$, an $(A,A_0)$-Morita
	bimodule $U$ gives an $(eAe,fA_0f)$-Morita bimodule $eUf$ by composing
	the Morita equivalences associated with $e$ and $f$; see
	\cite[Theorem~2.8.7 and Corollary~2.8.10]{Linckelmann2018a}.
	
	Suppose that a group $Z$ acts through homomorphisms
	$Z\longrightarrow Z(A)^\times$ and $Z\longrightarrow Z(A_0)^\times$.
	We say that $U$ is balanced over $\cO Z$ if $zu=uz$ for every $z\in Z$ and
	$u\in U$. For $\cO$-algebra automorphisms $a$ of $A$ and $a_0$ of $A_0$, an
	$\cO$-linear map $T:U\longrightarrow U$ is semilinear with respect to
	$(a,a_0)$ if
	\[
	T(xuy)=a(x)T(u)a_0(y)\qquad(x\in A,\ u\in U,\ y\in A_0).
	\]
	
	Let $\mu_{p'}(\cO)$ be the group of roots of unity in $\cO$ of order prime
	to $p$. Every scalar group $\mu$ below is a finite subgroup of
	$\mu_{p'}(\cO)$. We regard it as a trivial module for any group acting on
	it. Thus a normalized cocycle $\beta\in Z^2(\Lambda,\mu)$ satisfies
	\[
	\beta(x,y)\beta(xy,z)=\beta(y,z)\beta(x,yz),\qquad
	\beta(1,x)=\beta(x,1)=1.
	\]
	
	Let $b$ be a block of $\cO G$ with defect group $D$, and put
	\[
	N=N_G(D),\qquad
	\mathcal A_D=\Aut(G)_{D,b},\qquad
	I_D=\{c_n:n\in N\},\qquad
	\Omega_D=\mathcal A_D/I_D.
	\]
	Brauer correspondence implies that $\mathcal A_D$ also stabilizes $b_0$.
	Moreover, $I_D\cong N/Z(G)$. Defect-group conjugacy gives
	$\Aut(G)_b=\Inn(G)\mathcal A_D$, and
	$\Inn(G)\cap\mathcal A_D=I_D$. Hence
	\begin{equation}\label{eq:outer-defect-stabilizer}
		\Omega_D\cong\Aut(G)_b/\Inn(G).
	\end{equation}
	
	\begin{definition}\label{def:datum}
		A scalar-equivariant Morita datum for $(G,b,D)$ consists of an
		$(\cO Gb,\cO Nb_0)$-Morita bimodule $U$ balanced over $\cO Z(G)$, a finite
		subgroup $\mu\leq\mu_{p'}(\cO)$, a normalized cocycle
		$\beta\in Z^2(\Omega_D,\mu)$, and invertible maps
		$T_a:U\longrightarrow U$, for $a\in\mathcal A_D$, semilinear for the
		induced actions on the two block algebras, such that
		\begin{equation}\label{eq:scalar-datum}
			T_1=\id_U,\qquad
			T_{c_n}(u)=nun^{-1},\qquad
			T_aT_{a'}=\beta(aI_D,a'I_D)T_{aa'}
		\end{equation}
		for $a,a'\in\mathcal A_D$, $n\in N$, and $u\in U$.
	\end{definition}
	
	Central balance makes the prescribed map $u\mapsto nun^{-1}$ depend only
	on $c_n$. When $\beta=1$, the maps $T_a$ form a semilinear action of
	$\mathcal A_D$.
	
	\subsection{Tensor powers}
	
	\begin{lemma}[Tensor powers]\label{lem:tensor-datum}
		Let $X$ be a quasisimple group, and let $c$ be a block of $\cO X$
		with defect group $D$. Suppose that $(U,\mu,\beta,\{T_a\})$ is a
		scalar-equivariant Morita datum for $(X,c,D)$.
		For every positive integer $r$, put $Y=X^r$ and $D_Y=D^r$.
		Then $(Y,c^{\otimes r},D_Y)$ has a scalar-equivariant Morita datum
		with bimodule $U^{\otimes r}$ and scalar group $\mu$.
		The right block algebra is
		$\cO N_Y(D_Y)c_0^{\otimes r}$, where $c_0$ is the Brauer
		correspondent of $c$ in $N_X(D)$.
	\end{lemma}
	
	\begin{proof}
		Put $N=N_X(D)$,
		$\mathcal A=\Aut(X)_{D,c}$, $I=\{c_n:n\in N\}$, and
		$\Omega=\mathcal A/I$. The tensor product $U_Y=U^{\otimes r}$
		is a Morita bimodule between
		\[
		\cO Yc^{\otimes r}\qquad\text{and}\qquad
		\cO N^r c_0^{\otimes r},
		\]
		and is balanced over $\cO Z(Y)$. Moreover,
		$N_Y(D_Y)=N^r$, and $c_0^{\otimes r}$ is the Brauer
		correspondent of $c^{\otimes r}$.
		Since the $r$ coordinate subgroups of $Y$ are exactly its components,
		every automorphism of $Y$ permutes them. Hence
		\[
		\mathcal A_Y:=\Aut(Y)_{D_Y,c^{\otimes r}}
		=\mathcal A\wr S_r,\qquad
		I_Y:=\{c_y:y\in N_Y(D_Y)\}=I^r,
		\]
		and $\Omega_Y:=\mathcal A_Y/I_Y\cong\Omega\wr S_r$.
		For $q=((a_i),\sigma)\in\mathcal A_Y$, use the action
		$q((x_i))=(a_i(x_{\sigma^{-1}(i)}))$ on $Y$, and define
		\[
		T^Y_q(u_1\otimes\cdots\otimes u_r)
		=\bigotimes_{i=1}^r T_{a_i}(u_{\sigma^{-1}(i)}).
		\]
		These maps are invertible and semilinear for the actions of $q$
		on the two block algebras.
		For $q'=((b_i),\tau)$, set
		\begin{equation}\label{eq:wreath-multiplier}
			\beta^{\mathrm{wr}}(qI_Y,q'I_Y)
			=\prod_{i=1}^r\beta(a_iI,b_{\sigma^{-1}(i)}I).
		\end{equation}
		This defines a normalized cocycle in $Z^2(\Omega_Y,\mu)$:
		the cocycle identity follows by multiplying those for $\beta$ in
		the tensor factors and reindexing the factors permuted by $\sigma$.
		The identities \eqref{eq:scalar-datum} give
		\[
		T^Y_1=\id_{U_Y},\qquad
		T^Y_{c_y}(v)=yvy^{-1}\quad(y\in N_Y(D_Y),\ v\in U_Y),
		\]
		and
		\[
		T^Y_qT^Y_{q'}
		=\beta^{\mathrm{wr}}(qI_Y,q'I_Y)T^Y_{qq'}.
		\]
		Thus $(U_Y,\mu,\beta^{\mathrm{wr}},\{T^Y_q\})$ satisfies
		Definition~\ref{def:datum} for $(Y,c^{\otimes r},D_Y)$.
	\end{proof}
	
	\subsection{Extension from an orbit of components}
	
	We use the following standard extension criterion for graded
	Morita equivalences \cite[Theorem~5.1.2]{Marcus1999}.
	Let $C$ and $C_0$ be strongly $\Gamma$-graded algebras with
	identity components $A$ and $A_0$, respectively, and put
	\[
	\Delta(C,C_0)=\bigoplus_{g\in\Gamma}
	C_g\otimes((C_0)_{g^{-1}})^{\op}
	\subseteq C\otimes C_0^{\op}.
	\]
	If the $A\otimes A_0^{\op}$-action on an $(A,A_0)$-Morita
	bimodule $V$ extends to $\Delta(C,C_0)$, then
	$C\otimes_A V$ induces a graded Morita equivalence
	between $C$ and $C_0$.
	We also use the fact that the identity component of a graded
	Morita bimodule between strongly graded algebras induces a
	Morita equivalence between their identity components
	\cite[p.~74]{MarcusMinuta2020}.
	
	\begin{proposition}[Extension from a component orbit]\label{prop:quotient-induction}
		Let $b$ be a block of $\cO G$, and let $X$ be a component of $G$,
		that is, a subnormal quasisimple subgroup.
		Let $X_1=X,\ldots,X_r$ be the distinct $G$-conjugates of $X$,
		and put $L=X_1\cdots X_r$. Let $c$ be a $G$-stable block of
		$\cO L$ covered by $b$, let $D$ be a defect group of $c$, and
		let $f$ be the block of $\cO X$ covered by $c$.
		Suppose that $(X,f,D\cap X)$ admits a scalar-equivariant Morita
		datum with scalar group $\mu$, and put $M=N_G(D)$.
		Then there are a central extension
		\[
		1\longrightarrow\mu\longrightarrow\widehat M
		\xrightarrow{\rho}M\longrightarrow1
		\]
		and a block $b'$ of $\cO\widehat M$ such that
		\[
		\cO Gb\sim_{\Mor}\cO\widehat M b',\qquad l(b')=l(b).
		\]
		Moreover, the extension can be chosen to split over $Z(G)$, with
		\[
		\rho^{-1}(Z(G))=\mu\times Z(G)\leq Z(\widehat M).
		\]
	\end{proposition}
	
	\begin{proof}
		Put $N=N_L(D)$, and let $c_0$ be the Brauer correspondent of $c$
		in $N$. Since $c$ is $G$-stable, defect-group conjugacy gives
		$G=LM$ and $M/N\cong G/L$.
		\proofstep{Step 1. The tensor product bimodule.}
		For each $i$, let $f_i$ be the block of $\cO X_i$ covered by $c$,
		and put $D_i=D\cap X_i$. The block theory of central products
		gives $D=D_1\cdots D_r$, with $D_i$ a defect group of $f_i$
		\cite[Lemma~7.5(i)]{Sambale2014}.
		The group $L$ normalizes every $X_i$, so $M$ is transitive on
		$\{X_1,\ldots,X_r\}$. Choose $m_i\in M$ with $m_1=1$ and
		$m_iXm_i^{-1}=X_i$. The maps
		\[
		\iota_i:X\longrightarrow X_i,\qquad x\longmapsto m_ixm_i^{-1}
		\]
		identify $(X,f,D_1)$ with $(X_i,f_i,D_i)$.
		Choose a scalar-equivariant Morita datum
		$(U,\mu,\beta,\{T_a\})$ for $(X,f,D_1)$.
		Let $f_0$ be the Brauer correspondent of $f$ in $N_X(D_1)$, and put
		\[
		Y=X^r,\qquad D_Y=D_1^r,\qquad
		f_Y=f^{\otimes r},\qquad f_{Y,0}=f_0^{\otimes r}.
		\]
		By Lemma~\ref{lem:tensor-datum}, $U_Y=U^{\otimes r}$ is the
		bimodule of a scalar-equivariant Morita datum
		$(U_Y,\mu,\beta^{\mathrm{wr}},\{T^Y_q\})$ for $(Y,f_Y,D_Y)$.
		In particular, $U_Y$ is an
		$(\cO Yf_Y,\cO N_Y(D_Y)f_{Y,0})$-Morita bimodule balanced over
		$\cO Z(Y)$.
		Multiplication defines a surjection
		\[
		\pi:Y\longrightarrow L,\qquad
		(x_1,\ldots,x_r)\longmapsto\iota_1(x_1)\cdots\iota_r(x_r).
		\]
		Write $K_\pi=\ker(\pi)\leq Z(Y)$, to distinguish this kernel
		from the coefficient field $K$. The central-product description
		of $c$ shows that $f_Y$ dominates $c$ and $\pi(D_Y)=D$.
		We also use $\pi$ for the induced maps of group algebras.
		
		\proofstep{Step 2. Descent to the central quotient.}
		We now descend the tensor product bimodule through $\pi$.
		Write $K_\pi=K_{\pi,p}\times K_{\pi,p'}$ for its Sylow
		$p$-subgroup and Hall $p'$-subgroup. Every central $p$-subgroup
		is contained in every defect group, so $K_{\pi,p}\leq D_Y$.
		Thus $\pi^{-1}(D)=D_Y\times K_{\pi,p'}$ has unique Sylow
		$p$-subgroup $D_Y$. Hence
		\[
		\pi^{-1}(N_L(D))=N_Y(\pi^{-1}(D))=N_Y(D_Y),
		\]
		and consequently
		\begin{equation}\label{eq:product-normalizer}
			N_Y(D_Y)/K_\pi\cong N_L(D)=N.
		\end{equation}
		Put
		\[
		\begin{aligned}
			A_Y&=\cO Yf_Y,& A_{Y,0}&=\cO N_Y(D_Y)f_{Y,0},\\
			A&=\cO Lc,& A_0&=\cO Nc_0.
		\end{aligned}
		\]
		Compatibility of Brauer correspondence with successive central $p'$- and $p$-quotients \cite[Theorems~9.9(c) and~9.10]{Navarro1998} gives
		\[
		\pi(f_Y)=c,\qquad \pi(f_{Y,0})=c_0.
		\]
		For $z\in K_{\pi,p'}$ and $u\in U_Y$, we have
		\[
		zf_Y=f_Y,\qquad zf_{Y,0}=f_{Y,0},\qquad zu=u=uz.
		\]
		Thus the $p'$-part of the kernel already acts trivially.
		Define the ideals
		\[
		I=\sum_{z\in K_\pi}(z-1)A_Y,\qquad
		I_0=\sum_{z\in K_\pi}(z-1)A_{Y,0}.
		\]
		The restrictions of $\pi$ therefore induce algebra isomorphisms
		\begin{equation}\label{eq:central-quotient-algebras}
			A_Y/I\xrightarrow{\sim}A,\qquad
			A_{Y,0}/I_0\xrightarrow{\sim}A_0.
		\end{equation}
		Central balance gives the explicit equality
		\[
		IU_Y=\sum_{z\in K_\pi}(z-1)U_Y
		=\sum_{z\in K_\pi}U_Y(z-1)=U_Y I_0.
		\]
		The inverse Morita bimodule $U_Y^\vee$ is balanced as well, so
		$I_0 U_Y^\vee=U_Y^\vee I$. We can consequently define directly
		\begin{equation}\label{eq:central-quotient-bimodule}
			V=U_Y/IU_Y
			\cong A\otimes_{A_Y}U_Y
			\cong U_Y\otimes_{A_{Y,0}}A_0.
		\end{equation}
		Writing $\bar u=u+IU_Y$, the induced bimodule actions are
		\begin{equation}\label{eq:central-quotient-action}
			\pi(a)\bar u\,\pi(a_0)=\overline{au a_0}
			\qquad(a\in A_Y,\ a_0\in A_{Y,0},\ u\in U_Y).
		\end{equation}
		These actions are well defined because $IU_Y=U_Y I_0$.
		Set $V^\vee=U_Y^\vee/I_0 U_Y^\vee$. The identities
		\eqref{eq:morita-identities} descend to
		\[
		\begin{aligned}
			V\otimes_{A_0} V^\vee
			&\cong A\otimes_{A_Y}
			(U_Y\otimes_{A_{Y,0}}U_Y^\vee)\cong A,\\
			V^\vee\otimes_A V
			&\cong A_0\otimes_{A_{Y,0}}
			(U_Y^\vee\otimes_{A_Y}U_Y)\cong A_0.
		\end{aligned}
		\]
		Hence $V$ is an $(A,A_0)$-Morita bimodule; compare
		\cite[Remark~2.7]{KessarKoshitaniLinckelmann2015}. In particular, it is finitely generated
		projective over $A$ and therefore $\cO$-free.
		Since $Y$ is perfect, $\pi^{-1}(Z(L))=Z(Y)$: if
		$\pi(y)\in Z(L)$, the commutator map from $Y$ to $K_\pi$ determined
		by $y$ is a homomorphism, so it is trivial. Thus central balance
		on $Y$ descends to $Z(L)$.
		
		\proofstep{Step 3. Semilinear maps and the scalar cocycle.}
		The conjugation action of $M$ lifts to $Y$ as follows.
		Write $\mathcal A_Y=\Aut(Y)_{D_Y,f_Y}$ and
		$I_Y=\{c_y:y\in N_Y(D_Y)\}$. For $m\in M$, let
		$\sigma_m\in S_r$ satisfy $mX_jm^{-1}=X_{\sigma_m(j)}$, and set
		\[
		a_i(m)=\iota_i^{-1}\circ c_m\circ\iota_{\sigma_m^{-1}(i)}
		\in\Aut(X)_{D_1,f}.
		\]
		These maps define a homomorphism
		\[
		\nu:M\longrightarrow\mathcal A_Y,\qquad
		m\longmapsto((a_i(m)),\sigma_m),
		\]
		for which
		\begin{equation}\label{eq:lifted-conjugation}
			\pi\circ\nu(m)=c_m\circ\pi\qquad(m\in M).
		\end{equation}
		For $z\in K_\pi$ and $u\in U_Y$, semilinearity gives
		\[
		T^Y_{\nu(m)}((z-1)u)
		=(\nu(m)(z)-1)T^Y_{\nu(m)}(u).
		\]
		Since $\nu(m)(K_\pi)=K_\pi$, the submodule $IU_Y$ is stable
		under $T^Y_{\nu(m)}$. These maps therefore induce invertible
		semilinear maps $t_m$ on $V$.
		For $n\in N$, choose $y\in N_Y(D_Y)$ with $\pi(y)=n$, using
		\eqref{eq:product-normalizer}. Since distinct components commute,
		$\nu(n)=c_y$. Formula~\eqref{eq:central-quotient-action} then gives
		\[
		t_n(\bar u)=\overline{T^Y_{c_y}(u)}
		=\overline{yuy^{-1}}=n\bar u n^{-1}.
		\]
		In particular, $\nu(N)\leq I_Y$. Set
		\begin{equation}\label{eq:quotient-relative-data}
			\begin{gathered}
				\gamma(mN,m'N)=\beta^{\mathrm{wr}}(\nu(m)I_Y,\nu(m')I_Y),\\
				t_m(\bar u)=\overline{T^Y_{\nu(m)}(u)}\qquad(u\in U_Y)
			\end{gathered}
		\end{equation}
		The first formula defines a normalized cocycle
		$\gamma\in Z^2(M/N,\mu)$. Since the maps $\nu(m)$ preserve
		the quotient relations, the induced maps $t_m$ are invertible
		and semilinear for conjugation by $m$ on the two block algebras.
		The identities in \eqref{eq:scalar-datum} and $\nu(n)=c_y$ give
		\begin{equation}\label{eq:relative-datum}
			t_1=\id_V,\qquad
			t_n(v)=nvn^{-1},\qquad
			t_mt_{m'}=\gamma(mN,m'N)t_{mm'}
		\end{equation}
		for $n\in N$, $m,m'\in M$, and $v\in V$.
		
		\proofstep{Step 4. A graded Morita equivalence.}
		The cocycle $\gamma$ defines the central extension
		\[
		\widehat M=\mu\times_\gamma M,\qquad
		(\zeta,m)(\xi,m')=(\zeta\xi\gamma(mN,m'N),mm'),
		\]
		with $\rho(\zeta,m)=m$, and embed $N$ by $n\mapsto(1,n)$.
		For the inclusion $\varepsilon:\mu\hookrightarrow\cO^\times$, put
		\[
		e=e_{\varepsilon^{-1}}=\frac1{|\mu|}
		\sum_{\zeta\in\mu}\zeta(\zeta,1).
		\]
		Set $\Gamma=M/N\cong G/L$ and
		\[
		C=\cO Gc,\qquad
		C_0=\cO\widehat M c_0e.
		\]
		Normalization of $\gamma$ gives
		\[
		(\zeta,m)(1,n)(\zeta,m)^{-1}=(1,mnm^{-1})
		\qquad(\zeta\in\mu,\ m\in M,\ n\in N).
		\]
		Brauer correspondence makes $c_0$ $M$-stable, so $c_0e$ is
		central in $\cO\widehat M$. Also $(\zeta,1)e=\zeta^{-1}e$. The inverse central character
		produces the inverse multiplier on $C_0$, and hence the multiplier
		$\gamma$ in the diagonal action below.
		The natural strong $\Gamma$-grading on $C$ has $C_1=A$ and
		$C_{mN}=A(mc)$. The quotient $\widehat M/(\mu N)\cong\Gamma$
		grades $C_0$, with $(C_0)_1=\cO(\mu N)c_0e$ identified with $A_0$ by
		$nc_0\mapsto(1,n)c_0e$. For $m\in M$, set
		\[
		v_m=(1,m)c_0e\in C_0^\times.
		\]
		Then $(C_0)_{mN}=A_0v_m$, so this grading is strong, and
		\begin{equation}\label{eq:v-product}
			v_mv_{m'}=\gamma(mN,m'N)^{-1}v_{mm'}.
		\end{equation}
		Let $\Delta=\Delta(C,C_0)$ and define
		\[
		d_m=mc\otimes(v_m^{-1})^{\op}\in\Delta_{mN}.
		\]
		We have $\Delta_1=A\otimes A_0^{\op}$ and
		$\Delta_{mN}=\Delta_1d_m$. Multiplication in the opposite algebra
		reverses the order of the second factors, so \eqref{eq:v-product} gives
		\begin{align}
			d_md_{m'}&=\gamma(mN,m'N)d_{mm'},\label{eq:d-product}\\
			d_m(a\otimes a_0^{\op})d_m^{-1}
			&=(mam^{-1})\otimes(ma_0m^{-1})^{\op}
			\qquad(a\in A,\ a_0\in A_0).\label{eq:d-conjugation}
		\end{align}
		Let $a\otimes a_0^{\op}$ act on $V$ by $v\mapsto av a_0$, and let
		$d_m$ act by $t_m$. The multiplication law in
		\eqref{eq:relative-datum} agrees with \eqref{eq:d-product}, and
		semilinearity gives \eqref{eq:d-conjugation}. For $n\in N$, the
		element $d_n=nc\otimes(n^{-1}c_0)^{\op}$ already lies in
		$\Delta_1$, and its original action $v\mapsto nvn^{-1}$ is
		precisely $t_n$. These assignments therefore agree on $\Delta_1$.
		Choosing a transversal $S$ for $M/N$ containing $1$, they define
		a unital $\Delta$-module structure on $V$ by
		\[
		(xd_s)v=x\,t_s(v)\qquad(x\in\Delta_1,\ s\in S,\ v\in V).
		\]
		The diagonal-extension criterion now shows that $C\otimes_A V$
		induces a $G/L$-graded Morita equivalence between $C$ and $C_0$.
		In particular,
		\begin{equation}\label{eq:scalar-induction}
			\cO Gc\sim_{\Mor}\cO\widehat M c_0e_{\varepsilon^{-1}}.
		\end{equation}
		
		\proofstep{Step 5. Restriction to blocks and central splitting.}
		Since $b$ covers the $G$-stable block $c$, we have $bc=b$.
		The equivalence therefore restricts to a Morita equivalence between
		$\cO Gb$ and $\cO\widehat M b'$ for a block $b'$ of
		$\cO\widehat M$. Reduction modulo the maximal ideal preserves
		this equivalence, so $l(b')=l(b)$.
		
		Finally, if $z\in Z(G)$, then $\nu(z)=1$ because $z$ centralizes
		every $X_i$. The formula for $\gamma$ and its normalization give
		\[
		\gamma(zN,mN)=\gamma(mN,zN)=1\qquad(m\in M).
		\]
		Consequently the multiplication on $\rho^{-1}(Z(G))$ is the
		direct-product multiplication on $\mu\times Z(G)$, and this
		subgroup is central in $\widehat M$.
	\end{proof}

	\section{Blocks of quasisimple groups}\label{sec:quasisimple}
	Throughout this section, $X$ is quasisimple, $c$ is a block of
	$\cO X$, and $\Aut(X)_c$ is transitive on $\IBr(c)$. We fix a defect
	group $D$ of $c$, put $N=N_X(D)$, and denote the Brauer correspondent
	of $c$ in $N$ by $c_0$. Write
	\[
	\ker(c)=\bigcap_{\xi\in\Irr(c)}\ker(\xi),\qquad
	\Out(X)_c=\Aut(X)_c/\Inn(X).
	\]
	We call $c$ faithful if $\ker(c)=1$.
	The proof of Proposition~\ref{prop:local} proceeds through a
	classification dichotomy for faithful blocks, a construction in each
	of its two cases, and removal of the block kernel.
	
	\subsection{The classification dichotomy}
	
	\begin{proposition}\label{prop:classification}
		Assume that $\ker(c)=1$. Then the defect group $D$ is abelian
		and the block $c$ is inertial.
		Moreover, at least one of the following holds.
		\begin{enumerate}
			\item[(i)] There are a finite group $Y$ containing $X$ as a normal
			subgroup and a nilpotent block $f$ of $\cO Y$ covering $c$ such that
			$Y/X$ is an abelian $p'$-group, $c$ is $Y$-stable, and
			$Z(X)\leq Z(Y)$. Moreover, there is a homomorphism
			\begin{equation}\label{eq:local-automorphism-lift}
				\ell:\Aut(X)_c\longrightarrow\Aut(Y),\qquad
				\ell(a)|_X=a\quad(a\in\Aut(X)_c),
			\end{equation}
			satisfying $\ell(c_x)(y)=xyx^{-1}$ for $x\in X$ and $y\in Y$.
			\item[(ii)] The order of $\Out(X)_c$ is prime to $p$.
		\end{enumerate}
	\end{proposition}
	\begin{proof}
		By \cite[Proposition~1.2]{JiangZhangZhou2026}, $D$ is abelian. Applying
		\cite[Proposition~4.1]{JiangZhangZhou2026} with $H=\widetilde H=X$,
		$\widetilde c=c$, and $\widetilde P=D$ then shows that $c$ is inertial.
		If $c$ is nilpotent, take
		$Y=X$, $f=c$, and $\ell$ to be the inclusion. We may therefore
		assume that $c$ is nonnilpotent and $p\mid|\Out(X)_c|$, and prove (i).
		
		Set $S=X/Z(X)$, $Z_p=O_p(Z(X))$, and $X^*=X/Z_p$.
		Since $O_p(X)=Z_p$ and $\ker(c)=1$, \cite[Corollary~6.12]{Navarro1998} gives
		$\bigcap_{\varphi\in\IBr(c)}\ker(\varphi)=Z_p$. By
		\cite[Theorem~9.10]{Navarro1998}, $c$ dominates a unique block
		$c^*$ of $X^*$, and inflation identifies $\IBr(c^*)$ with $\IBr(c)$.
		Thus $c^*$ is faithful. It is nonnilpotent by
		\cite[Proposition~2.7]{AnEaton2013}.
		Since $X$ is quasisimple and $Z_p\leq Z(X)$, we have
		\[
		Z(X^*)=Z(X)/Z_p,\qquad X^*/Z(X^*)\cong S.
		\]
		Moreover, $Z_p$ is characteristic in $X$, so the quotient maps
		induce the natural embeddings
		\begin{equation}\label{eq:local-central-p-quotient}
			\Out(X)_c\hookrightarrow\Out(X^*)_{c^*},\qquad
			\Out(X)\hookrightarrow\Out(S).
		\end{equation}
		These follow from \cite[Corollary~B.8]{Navarro2018}, applied to a common
		universal cover.
		
		If $S$ is sporadic, alternating or the Tits group, the classification
		in \cite[Proposition~2.6 and Theorem~2.8]{MalleNavarroSpath2018}, together
		with the nonnilpotency of $c^*$, gives $p\in\{3,5\}$.
		The corresponding groups $\Out(X^*)$ are $2$-groups \cite[proof of Lemma~4.3]{JiangZhangZhou2026}.
		If $X^*$ is an exceptional covering group, then, since $c^*$ is
		faithful and nonnilpotent, and $\Aut(X^*)_{c^*}$ is transitive on
		$\IBr(c^*)$,
		\cite[Proposition~2.11]{MalleNavarroSpath2018} gives $X^*\cong 2.G_2(4)$ and $p=3$.
		In this case $\Out(X^*)\cong C_2$ \cite[p.~97]{ConwayEtAl1985}.
		In either case, \eqref{eq:local-central-p-quotient} yields
		$p\nmid|\Out(X)_c|$, so (ii) holds, contrary to our assumption.
		The defining-characteristic case is excluded by
		\cite[Proposition~3.1]{MalleNavarroSpath2018}, since it forces $c^*$ to have central
		defect and hence to be nilpotent. Thus $S$ is of Lie type in
		characteristic $r\ne p$, and $X^*$ is not an exceptional covering
		group. Let $\widehat S$ be the universal covering group of $S$,
		and let $E$ be the kernel of its map onto the usual simply connected
		covering group. Then $E$ is an $r$-group \cite[p.~6 of the author manuscript]{FengEtAl2019}.
		Since $X^*$ is nonexceptional, $E$ is contained in the kernel of
		the composite homomorphism $\widehat S\longrightarrow X\longrightarrow X^*$.
		Hence the image $E_X$ of $E$ in $X$ lies in $Z_p$.
		But $E_X$ is an $r$-group and $Z_p$ is a $p$-group, so $E_X=1$.
		Thus the map $\widehat S\longrightarrow X$ factors through the
		simply connected covering group, and $X$ is nonexceptional as well.
		
		Suppose first that $S$ is a Suzuki--Ree group. Since $X$ is
		nonexceptional, $X=S$. By the proof of \cite[Lemma~4.4]{JiangZhangZhou2026},
		the nonnilpotency of $c$ leaves only $S\cong {}^2G_2(3^{2a+1})$
		with $p=2$. But $\Out(S)$ has odd order, contradicting
		$p\mid|\Out(X)_c|$. Thus the Suzuki--Ree families are excluded.
		
		We next exclude the untwisted families $B_2(2^f)$, $G_2(3^f)$
		and $F_4(2^f)$. Their simply connected algebraic groups have
		trivial center \cite[p.~3 of the arXiv version]{Johansson2022}, so again $X=S$ since $X$
		is nonexceptional. The identity embedding is therefore a regular
		embedding. Taking both finite groups equal to $X$ and both blocks
		equal to $c$ in \cite[Proposition~3.8]{JiangZhangZhou2026}, we obtain that $c$
		is nilpotent, contrary to our assumption.
		
		In the remaining cases, write $X=G/Z$ with $Z\leq Z(G)$,
		where $G=\mathbf G^F$ is perfect, $\mathbf G$ is a simple simply
		connected algebraic group in characteristic $r$, and $F$ is a
		Frobenius endomorphism, as in the proof of \cite[Proposition~4.2]{JiangZhangZhou2026}.
		Let $b$ be the block of $G$ dominating $c$. Inflation identifies
		$\IBr(c)$ with $\IBr(b)$. The center of $G$ has order prime to $r$,
		so $E$ is the Sylow $r$-subgroup of $Z(\widehat S)$ and is characteristic in $\widehat S$.
		Every automorphism of $X$ therefore lifts through $G$, by \cite[Corollary~B.8]{Navarro2018}.
		The lift is unique because $G$ is perfect and $Z$ is central.
		Thus $\Aut(G)_Z\cong\Aut(X)$ compatibly with inflation, and $\IBr(b)$ is a single
		$\Aut(G)_b$-orbit. Moreover, $b$ is nonnilpotent, since otherwise
		its central quotient $c$ would be nilpotent.
		
		Choose a regular embedding
		\[
		\iota:\mathbf G\hookrightarrow\widetilde{\mathbf G}
		\]
		with compatible field and graph actions as in
		\cite[pp.~873--874]{MalleSpath2016}. Identifying $\mathbf G$ with its image,
		we have $\widetilde{\mathbf G}=\mathbf G Z(\widetilde{\mathbf G})$
		and $Z(\widetilde{\mathbf G})$ is connected. Extend $F$ to
		$\widetilde{\mathbf G}$ and set $\widetilde G=\widetilde{\mathbf G}^{F}$.
		Choose a block $\widetilde b$ of $\widetilde G$ covering $b$ in
		the corresponding Lusztig series. Choose the semisimple $p'$-labels in the two dual groups
		compatibly, as required in \cite[Proposition~3.8]{JiangZhangZhou2026}.
		That proposition says that $b$ or $\widetilde b$ is nilpotent.
		Since $b$ is nonnilpotent, $\widetilde b$ is nilpotent.
		
		Let $\mathcal D\leq\Aut(G)$ be generated by the standard field
		and graph automorphisms commuting with $F$. Their compatible
		extensions act on $\widetilde G$, and the map
		\[
		\widetilde G\rtimes\mathcal D\longrightarrow\Aut(G),\qquad
		(g,d)\longmapsto\bigl(x\mapsto g\,d(x)g^{-1}\bigr)
		\]
		is surjective with kernel $Z(\widetilde G)\times\{1\}$
		\cite[pp.~873--874]{MalleSpath2016}. Put $\widetilde X=\widetilde G/Z$ and
		$\mathcal D_Z=\{d\in\mathcal D\mid d(Z)=Z\}$. Since
		$Z\leq Z(G)\leq Z(\widetilde G)$, the preimage of $\Aut(G)_Z$
		is $\widetilde G\rtimes\mathcal D_Z$. The identification
		$\Aut(G)_Z\cong\Aut(X)$ therefore gives an epimorphism
		\[
		\pi:\widetilde G\rtimes\mathcal D_Z\longrightarrow\Aut(X),
		\qquad \pi(g,d)(uZ)=g\,d(u)g^{-1}Z\quad(u\in G),
		\]
		with $\ker\pi=Z(\widetilde G)\times\{1\}$.
		Define a second homomorphism by
		\[
		\rho:\widetilde G\rtimes\mathcal D_Z\longrightarrow\Aut(\widetilde X),
		\qquad \rho(g,d)(yZ)=g\,d(y)g^{-1}Z\quad(y\in\widetilde G).
		\]
		For $z\in Z(\widetilde G)$, we have
		$\rho(z,1)(yZ)=zyz^{-1}Z=yZ$ for every $y\in\widetilde G$.
		Hence $\ker\pi\leq\ker\rho$. There is consequently a unique
		homomorphism
		\[
		\widetilde\ell:\Aut(X)\longrightarrow\Aut(\widetilde X),
		\qquad \widetilde\ell\bigl(\pi(g,d)\bigr)=\rho(g,d).
		\]
		The kernel inclusion makes this definition independent of the
		choice of $(g,d)$. Since $\rho(g,d)|_X=\pi(g,d)$, we have
		$\widetilde\ell(a)|_X=a$ for all $a\in\Aut(X)$.
		For $x=uZ\in X$ with $u\in G$, the element $(u,1)$ maps under
		$\pi$ to $c_x$, while $\rho(u,1)$ is conjugation by $x$ on
		$\widetilde X$. Thus $\widetilde\ell(c_x)$ is conjugation by $x$.
		
		Since $\ker(c)=1$, we have $\ker(b)=O_{p'}(Z)$. This subgroup
		acts trivially on the characters of $b$ and is central in
		$\widetilde G$. Restriction from a character of $\widetilde b$
		therefore lies over the same trivial central character, giving
		\[
		O_{p'}(Z)=\ker(b)\leq\ker(\widetilde b).
		\]
		Taking the quotient first by $O_{p'}(Z)$ and then by the central
		$p$-subgroup $Z/O_{p'}(Z)$ yields a nilpotent block of
		$\widetilde X=\widetilde G/Z$ covering $c$.
		The group $\widetilde X/X\cong\widetilde G/G$ is abelian.
		Let $Y'$ be the inverse image of $O_{p'}(\widetilde X/X)$.
		A block of $Y'$ covered by this nilpotent block
		and covering $c$ is nilpotent, since $\widetilde X/Y'$ is a
		$p$-group \cite[Lemma~2.2]{AnEaton2011}. Its Fong--Reynolds correspondent
		$f$ in $Y=I_{Y'}(c)$ is nilpotent and covers the $Y$-stable block $c$.
		Now $Y/X$ is an abelian $p'$-group and $Z(X)\leq Z(Y)$.
		The action of $\Aut(X)_c$ through $\widetilde\ell$ preserves $Y'$
		and $Y$; its restriction to $Y$ is the required homomorphism $\ell$,
		with $\ell(c_x)$ equal to conjugation by $x$. This proves (i).
	\end{proof}
	
	\subsection{Descent from a nilpotent covering block}
	
	\begin{lemma}[Nilpotent covering case]\label{lem:nilpotent-covering}
		Suppose that $Y,f,\ell$ satisfy the conditions in
		Proposition~\ref{prop:classification}\textup{(i)}.
		Then $(X,c,D)$ admits a scalar-equivariant Morita datum in the
		sense of Definition~\ref{def:datum}, with $\mu=\{1\}$ and $\beta=1$.
	\end{lemma}
	
	\begin{proof}
		Choose $Y,f,\ell$ as in Proposition~\ref{prop:classification}\textup{(i)},
		and put
		\[
		M=N_Y(D),\qquad N=N_X(D),\qquad \Gamma=Y/X.
		\]
		The Frattini argument gives $Y=XM$, and we identify $M/N$ with
		$\Gamma$ via $mN\mapsto mX$. Retain the Brauer correspondent
		$c_0$ of $c$ in $N$, and set
		\[
		\begin{aligned}
			A&=\cO Xc,& A_0&=\cO Nc_0,\\
			C&=\cO Yc,& C_0&=\cO Mc_0.
		\end{aligned}
		\]
		Here $c$ is $Y$-stable, and Brauer correspondence makes $c_0$
		$M$-stable. Consequently $C$ and $C_0$ have identities $c$ and
		$c_0$, respectively.
		
		We write $\alpha\in\Aut(X)_c$ for an automorphism of $X$ and
		$a=\ell(\alpha)$ for its lift to $Y$. Since $\ell(\alpha)|_X=\alpha$, the
		homomorphism $\ell$ is injective. Write
		\[
		\mathcal A_D=\Aut(X)_{D,c},\qquad
		\mathcal B=\ell(\Aut(X)_c),\qquad
		\mathcal B_D=\ell(\mathcal A_D).
		\]
		Thus $\mathcal B_D$ is the stabilizer of $D$ in $\mathcal B$;
		it preserves $M,N,c$, and $c_0$.
		
		Choose representatives $m_\gamma\in M$ for $\gamma\in\Gamma$,
		so that $m_\gamma X=\gamma$, with $m_1=1$.
		The quotient grading on $C$ is
		\[
		C=\bigoplus_{\gamma\in\Gamma}C_\gamma,\qquad
		C_\gamma=A(m_\gamma c)=\cO X\,m_\gamma c,\qquad C_1=A.
		\]
		The element $m_\gamma c$ is the image of $m_\gamma$ in $C$ under
		multiplication by $c$. Since $c$ commutes with $Y$, the element
		$m_\gamma c$ is a unit of $C$ with inverse $m_\gamma^{-1}c$.
		Similarly, $C_0$ is strongly $\Gamma$-graded with
		$(C_0)_\gamma=A_0(m_\gamma c_0)$ and $(C_0)_1=A_0$.
		
		Let $\mathcal F=\Hom(\Gamma,\cO^\times)$, with trivial character
		$1_{\mathcal F}$. For $\chi\in\mathcal F$,
		the character twist $\tau_\chi$ is the $\cO$-algebra automorphism
		of $\cO Y$ defined by
		\[
		\tau_\chi(y)=\chi(yX)y\qquad(y\in Y).
		\]
		It fixes $\cO X$ pointwise and restricts to $C$ and $C_0$.
		The group $\mathcal B$ acts on $\mathcal F$ by
		\[
		(a\cdot\chi)(yX)=\chi(a^{-1}(y)X).
		\]
		Since $a\tau_\chi a^{-1}=\tau_{a\cdot\chi}$, the rule
		$(\chi,a)\mapsto\tau_\chi\circ a$ defines an action of
		$\mathcal F\rtimes\mathcal B$ on $\cO Y$.
		In particular,
		\[
		H=\mathcal F\rtimes\mathcal B_D
		\]
		acts on both $C$ and $C_0$.
		
		\proofstep{Step 1. Covering blocks and their stabilizers.}
		Every block of $Y$ covering $c$ has $D$ as a defect group,
		because $Y/X$ has order prime to $p$.
		Character twists are compatible with the Brauer maps and commute
		with conjugation by $Y$, since
		$\chi(tyt^{-1}X)=\chi(yX)$ for $t,y\in Y$.
		They therefore transport Brauer pairs with unchanged stabilizers
		and preserve nilpotency.
		The Harris--Kn\"orr correspondence
		\cite[Theorem~9.28]{Navarro1998}, is consequently an
		$H$-equivariant bijection between the blocks of $Y$ covering $c$
		and those of $M$ covering $c_0$.
		By \cite[Lemma~2.2(a)]{KessarMalle2013}, $\mathcal F$ is transitive on the
		first set, hence on both sets.
		
		Let $f_0$ be the Brauer correspondent of the chosen nilpotent
		block $f$ in $M$, and put
		\[
		J=\operatorname{Stab}_{\mathcal F\rtimes\mathcal B}(f),
		\qquad J_D=J\cap H.
		\]
		Equivariance shows that $f$ and $f_0$ have the same stabilizer
		$J_D$ in $H$.
		
		\proofstep{Step 2. A semilinear action on one covering block.}
		To construct an equivariant Morita bimodule for $f,f_0$,
		we realize the character twists by group automorphisms
		and apply the nilpotent-extension theorem.
		Let $Z_1\leq\cO^\times$ be the group of $\exp(\Gamma)$-th roots
		of unity, and let $\rho:Z_1\hookrightarrow\cO^\times$ be the inclusion.
		Our splitting assumptions give $|Z_1|=\exp(\Gamma)$, which is
		prime to $p$, and every $\chi\in\mathcal F$ takes values in $Z_1$.
		On $\widetilde Y=Y\times Z_1$, define
		\[
		\widetilde\tau_\chi(y,z)=(y,\chi(yX)z),\qquad
		\widetilde a(y,z)=(a(y),z)\quad(a\in\mathcal B).
		\]
		These maps define an action of $\mathcal F\rtimes\mathcal B$
		on $\widetilde Y$, fixing $Z_1$ pointwise.
		The central idempotent
		\[
		e_\rho=\frac1{|Z_1|}\sum_{z\in Z_1}\rho(z)^{-1}(1,z)
		\]
		is fixed, and the algebra isomorphism
		\[
		\cO\widetilde Y e_\rho\xrightarrow{\sim}\cO Y,\qquad
		(y,z)e_\rho\longmapsto\rho(z)y
		\]
		is equivariant for this action and the character-twist action.
		We use this isomorphism to identify the corresponding block algebras.
		
		Set $T=\widetilde Y\rtimes J$, using the group action just defined.
		We use the standard embeddings $\widetilde Y,J\leq T$ and identify
		$D$ with $D\times1\leq\widetilde Y$.
		The block $fe_\rho$ is $T$-stable and nilpotent, with defect group
		$D$, and its Brauer correspondent in
		$N_{\widetilde Y}(D)=M\times Z_1$ is $f_0e_\rho$.
		Defect-group conjugacy gives $T=\widetilde Y N_T(D)$, so the
		quotient on the normalizer side is naturally identified with
		$T/\widetilde Y$.
		The nilpotent-extension theorem \cite[Theorem~3.1]{Rossi2025}, following
		\cite[Theorem~3.14 and Corollary~3.15]{PuigZhou2012}, gives a
		$T/\widetilde Y$-graded basic Morita equivalence between
		\[
		\cO T(fe_\rho)\qquad\text{and}\qquad
		\cO N_T(D)(f_0e_\rho).
		\]
		Choose a graded Morita bimodule $W_f$ for this equivalence and
		put $V_f=(W_f)_1$, its identity component for the
		$T/\widetilde Y$-grading. Thus $V_f$ is a
		$(\cO Yf,\cO Mf_0)$-Morita bimodule.
		For $g\in N_T(D)$, define
		\begin{equation}\label{eq:local-graded-conjugation}
			t_g(v)=gvg^{-1}\qquad(v\in V_f),
		\end{equation}
		using the left and right actions on $W_f$, with the block
		idempotents suppressed.
		If $g$ has degree $g\widetilde Y$, then $gvg^{-1}$ has degree $1$,
		so $t_g$ preserves $V_f$.
		The bimodule identities give
		\[
		t_1=\id_{V_f},\qquad t_gt_{g'}=t_{gg'},\qquad
		t_g(rvs)=(grg^{-1})t_g(v)(gsg^{-1}),
		\]
		where $r\in\cO\widetilde Y(fe_\rho)$ and
		$s\in\cO(M\times Z_1)(f_0e_\rho)$, using the algebra
		isomorphisms above.
		Every element of $J_D$ stabilizes $D$, so $J_D\leq N_T(D)$
		acts semilinearly on $V_f$;
		we denote this action by $t_j$, for $j\in J_D$.
		
		\proofstep{Step 3. Normalization on inner automorphisms.}
		We must also prescribe the action of inner automorphisms.
		For $x\in X$, set
		\[
		\widehat c_x=\ell(c_x)\in\mathcal B,\qquad
		\iota_x=(1_{\mathcal F},\widehat c_x)\in\mathcal F\rtimes\mathcal B.
		\]
		By the property of $\ell$, $\widehat c_x$ is conjugation by $x$
		on $Y$. Thus $\iota_x$ fixes $f$, so $\iota_x\in J$, and it acts
		on $\widetilde Y$ by conjugation by $(x,1)$.
		Consequently
		\[
		\delta_x=((x^{-1},1),\iota_x)\in T
		\]
		lies in $C_T(\widetilde Y)\leq N_T(D)$.
		Semidirect-product multiplication gives
		$\delta_x\delta_{x'}=\delta_{xx'}$.
		It follows that $x\mapsto t_{\delta_x}$ is a homomorphism
		\[
		X\longrightarrow
		\Aut_{\cO Yf\text{-}\cO Mf_0}(V_f)
		\cong Z(\cO Yf)^\times.
		\]
		The target is abelian and $X$ is perfect, so $t_{\delta_x}=\id_{V_f}$.
		For $n\in N$, we have $\iota_n\in J_D$ and
		$\iota_n=(n,1)\delta_n$ in $N_T(D)$. Hence
		\begin{equation}\label{eq:local-cover-inner}
			t_{\iota_n}(v)=nvn^{-1}\qquad(n\in N,\ v\in V_f).
		\end{equation}
		\proofstep{Step 4. Induction over all covering blocks.}
		Regard $V_f$ as a $(C,C_0)$-bimodule via the block projections
		$r\mapsto rf$ and $s\mapsto sf_0$.
		The maps $t_j$ remain semilinear for these actions.
		Form
		\[
		V=\cO H\otimes_{\cO J_D}V_f.
		\]
		For $r\in C$, $s\in C_0$, $h,k\in H$, and $v\in V_f$, set
		\[
		\begin{aligned}
			r(h\otimes v)s&=h\otimes h^{-1}(r)v h^{-1}(s),\\
			S_k(h\otimes v)&=kh\otimes v.
		\end{aligned}
		\]
		The bimodule formula respects the tensor relation
		$hj\otimes v=h\otimes t_j(v)$, since
		\[
		\begin{aligned}
			r(hj\otimes v)s
			&=h\otimes t_j\bigl(j^{-1}h^{-1}(r)v j^{-1}h^{-1}(s)\bigr)\\
			&=h\otimes h^{-1}(r)t_j(v)h^{-1}(s)
			=r(h\otimes t_j(v))s.
		\end{aligned}
		\]
		The maps $S_k$ give the usual induced $\cO$-linear $H$-action:
		$S_1=\id_V$ and $S_kS_{k'}=S_{kk'}$.
		Moreover,
		\[
		\begin{aligned}
			S_k\bigl(r(h\otimes v)s\bigr)
			&=kh\otimes h^{-1}(r)v h^{-1}(s)\\
			&=kh\otimes (kh)^{-1}(k(r))v (kh)^{-1}(k(s))\\
			&=k(r)S_k(h\otimes v)k(s),
		\end{aligned}
		\]
		where the second equality uses $(kh)^{-1}\circ k=h^{-1}$.
		Thus $S_k(rws)=k(r)S_k(w)k(s)$ for $w\in V$:
		the action is semilinear, so $V$ is an $H$-equivariant
		$(C,C_0)$-bimodule.
		
		Let $\mathcal R\subseteq H$ represent the left cosets of $J_D$
		in $H$. The right $\cO J_D$-module decomposition
		$\cO H=\bigoplus_{h\in\mathcal R}h\,\cO J_D$ gives
		\[
		V=\bigoplus_{h\in\mathcal R}(h\otimes V_f).
		\]
		Each summand transports the original Morita equivalence to the
		corresponding block algebras
		$\cO Y\,h(f)$ and $\cO M\,h(f_0)$.
		Transitivity and the common stabilizer $J_D$ ensure that each
		corresponding pair occurs exactly once.
		Since $C$ and $C_0$ are the sums of their covering block algebras,
		$V$ is therefore a $(C,C_0)$-Morita bimodule.
		
		To check the inner normalization, write $h=(\chi,a)\in H$.
		Character twists fix $N$ pointwise, so $h^{-1}(n)=a^{-1}(n)$
		for $n\in N$. Also, $\widehat c_n$ fixes $\mathcal F$, and
		$\widehat c_n a=a\widehat c_{a^{-1}(n)}$.
		Multiplication in $H$ therefore gives
		$\iota_n h=h\iota_{a^{-1}(n)}$, with $\iota_{a^{-1}(n)}\in J_D$.
		Using the tensor relation and \eqref{eq:local-cover-inner}, we get
		\[
		\begin{aligned}
			S_{\iota_n}(h\otimes v)
			&=h\otimes t_{\iota_{a^{-1}(n)}}(v)\\
			&=h\otimes a^{-1}(n)v a^{-1}(n^{-1})\\
			&=n(h\otimes v)n^{-1}.
		\end{aligned}
		\]
		Hence $S_{\iota_n}(w)=nwn^{-1}$ for all $n\in N$ and $w\in V$.
		
		\proofstep{Step 5. Recovery of the identity components.}
		To pass from $C,C_0$ to $A,A_0$, we extend this equivalence
		to the skew group algebras and then apply averaging idempotents.
		For $\chi\in\mathcal F$ and $a\in\mathcal B_D$, abbreviate
		\[
		S_\chi=S_{(\chi,1)},\qquad S_a=S_{(1_{\mathcal F},a)}.
		\]
		Write $C\#_\tau\cO\mathcal F$ for the skew group algebra,
		with $r\#_\tau\chi=r\otimes\chi$ in its underlying module
		$C\otimes_\cO\cO\mathcal F$.
		For $r,r'\in C$ and $\chi,\psi\in\mathcal F$, multiplication is
		\[
		(r\#_\tau\chi)(r'\#_\tau\psi)
		=r\tau_\chi(r')\#_\tau\chi\psi.
		\]
		We embed $C$ by $r\mapsto r\#_\tau1_{\mathcal F}$ and suppress
		this embedding in products. Use the same conventions for $C_0$.
		The right scalar extension
		\[
		W=V\otimes_{C_0}(C_0\#_\tau\cO\mathcal F)
		\cong V\otimes_\cO\cO\mathcal F
		\]
		has bimodule actions
		\[
		\begin{aligned}
			(r\#_\tau\chi)(w\otimes\psi)
			&=rS_\chi(w)\otimes\chi\psi,\\
			(w\otimes\psi)(s\#_\tau\chi)
			&=w\tau_\psi(s)\otimes\psi\chi,
		\end{aligned}
		\]
		for $r\in C$, $s\in C_0$, $w\in V$, and $\chi,\psi\in\mathcal F$.
		Semilinearity of $S_\chi$ and its group law verify the bimodule
		identities. Equivariant Morita theory
		\cite[Theorem~5.1.2]{Marcus1999} shows that $W$ is a
		$(C\#_\tau\cO\mathcal F,C_0\#_\tau\cO\mathcal F)$-Morita bimodule.
		
		Extend the $\mathcal B_D$-action to $C\#_\tau\cO\mathcal F$ by
		$a(r\#_\tau\chi)=a(r)\#_\tau(a\cdot\chi)$, and similarly on
		$C_0\#_\tau\cO\mathcal F$. On $W$, define
		\[
		\widetilde S_a(w\otimes\chi)=S_a(w)\otimes(a\cdot\chi).
		\]
		The relations $S_aS_\chi=S_{a\cdot\chi}S_a$ and
		$a\tau_\chi=\tau_{a\cdot\chi}a$ give semilinearity, and the group
		law on $V$ gives
		$\widetilde S_a\widetilde S_b=\widetilde S_{ab}$.
		Thus $\{\widetilde S_a:a\in\mathcal B_D\}$ is a semilinear
		action on $W$.
		
		Since $\Gamma$ is abelian of order prime to $p$, our splitting
		assumptions give $|\mathcal F|=|\Gamma|\in\cO^\times$.
		Define the averaging idempotent
		\[
		p_1=\frac1{|\mathcal F|}\sum_{\chi\in\mathcal F}
		(1_C\#_\tau\chi).
		\]
		The group law on $\mathcal F$ gives $p_1^2=p_1$.
		For $r=\sum_{\gamma\in\Gamma}r_\gamma\in C$, with
		$r_\gamma\in C_\gamma$, character orthogonality gives
		$|\mathcal F|^{-1}\sum_{\chi\in\mathcal F}\tau_\chi(r)=r_1\in A$.
		Since $(1_C\#_\tau\psi)p_1=p_1$ for every $\psi\in\mathcal F$,
		\[
		p_1(r\#_\tau\psi)p_1
		=\left(\frac1{|\mathcal F|}
		\sum_{\chi\in\mathcal F}\tau_\chi(r)\right)p_1=r_1p_1.
		\]
		Consequently
		\begin{equation}\label{eq:local-full-corners}
			p_1(C\#_\tau\cO\mathcal F)p_1=Ap_1\cong A.
		\end{equation}
		The isomorphism is $r\mapsto rp_1$: it is multiplicative because
		$p_1$ centralizes $A$, and injective because the coefficient of
		$1_{\mathcal F}$ in $rp_1$ is $r/|\mathcal F|$.
		
		For the homogeneous units $m_\gamma c$ chosen above,
		$\tau_\chi(m_\gamma c)=\chi(\gamma)m_\gamma c$.
		Conjugation by $m_\gamma c$ therefore multiplies the $\chi$-term
		in $p_1$ by $\chi(\gamma)^{-1}$, and character orthogonality gives
		\[
		\begin{aligned}
			\sum_{\gamma\in\Gamma}(m_\gamma c)p_1(m_\gamma c)^{-1}
			&=\frac1{|\mathcal F|}\sum_{\chi\in\mathcal F}
			\left(\sum_{\gamma\in\Gamma}\chi(\gamma)^{-1}\right)
			(1_C\#_\tau\chi)\\
			&=1.
		\end{aligned}
		\]
		Thus $p_1$ is full.
		Define $p_1^{(0)}$ by the same averaging formula in
		$C_0\#_\tau\cO\mathcal F$, with $1_{C_0}$ in place of $1_C$.
		The same argument shows that $p_1^{(0)}$ is full and identifies
		$p_1^{(0)}(C_0\#_\tau\cO\mathcal F)p_1^{(0)}$ with $A_0$.
		Composing the equivalence induced by $W$ with the Morita
		equivalences defined by $p_1$ and $p_1^{(0)}$, as in
		\cite[Theorem~2.8.7 and Corollary~2.8.10]{Linckelmann2018a}, gives an
		$(A,A_0)$-Morita bimodule
		\[
		U=p_1Wp_1^{(0)}.
		\]
		In particular, $U$ is finitely generated projective over $A$,
		hence $\cO$-free.
		
		\proofstep{Step 6. The resulting datum.}
		Each $a\in\mathcal B_D$ permutes $\mathcal F$, so it fixes $p_1$
		and, similarly, $p_1^{(0)}$. Hence $\widetilde S_a(U)=U$.
		For $\alpha\in\mathcal A_D$, define
		\[
		T_\alpha=\left.\widetilde S_{\ell(\alpha)}\right|_U.
		\]
		These maps are semilinear for the actions induced by $\alpha$ on
		$A,A_0$, and satisfy
		\[
		T_1=\id_U,\qquad T_\alpha T_{\alpha'}=T_{\alpha\alpha'}
		\qquad(\alpha,\alpha'\in\mathcal A_D).
		\]
		For $n\in N$, the lift $\widehat c_n=\ell(c_n)$ fixes $\Gamma$
		and hence $\mathcal F$. Since $S_{\widehat c_n}=S_{\iota_n}$ on $V$
		and character twists fix $N$, the identity
		$S_{\iota_n}(w)=nwn^{-1}$ and the definition of $\widetilde S_a$ give
		$\widetilde S_{\widehat c_n}(\omega)=n\omega n^{-1}$ for $\omega\in W$.
		Restricting to $U$ yields
		\[
		T_{c_n}(u)=nun^{-1}\qquad(n\in N,\ u\in U).
		\]
		For $z\in Z(X)\leq N$, we have $c_z=1$, so the last identity gives
		$zu=uz$. Thus $U$ is balanced over $\cO Z(X)$.
		All conditions of Definition~\ref{def:datum} hold with
		$\mu=\{1\}$ and $\beta=1$.
	\end{proof}
	
	\subsection{An outer stabilizer of order prime to \texorpdfstring{$p$}{p}}
	
	\begin{lemma}[Outer stabilizer of order prime to $p$]\label{lem:outer-pprime}
		Under the hypotheses of Proposition~\ref{prop:classification},
		assume that \textup{(ii)} holds. Then $(X,c,D)$ admits a
		scalar-equivariant Morita datum in the sense of
		Definition~\ref{def:datum}.
	\end{lemma}
	\begin{proof}
		By Proposition~\ref{prop:classification}, $D$ is abelian and $c$
		is inertial. Let $N=N_X(D)$ and let $c_0$ be the Brauer
		correspondent of $c$ in $N$. Put
		\[
		A=\cO Xc,\qquad A_0=\cO Nc_0,\qquad
		\mathcal A=\Aut(X)_c,\qquad I=\Inn(X).
		\]
		Retain $\mathcal A_D=\Aut(X)_{D,c}$, $I_D=\{c_n:n\in N\}$,
		and $\Omega_D=\mathcal A_D/I_D$.
		
		For $P\in\Syl_p(I)$, the Frattini argument gives
		$\mathcal A=I N_{\mathcal A}(P)$.
		Since $|\mathcal A:I|$ is prime to $p$, $P$ is a normal Hall
		subgroup of $N_{\mathcal A}(P)$.
		Schur--Zassenhaus gives a complement $H$ to $P$ in
		$N_{\mathcal A}(P)$, of order prime to $p$, with $\mathcal A=IH$.
		This does not assert that $H\cap I=1$; that intersection is
		accounted for by the kernel $\mathcal K$ below. Set
		\[
		G=X\rtimes H,\qquad M=N_G(D),\qquad \bar{G}=G/X,
		\]
		and let $\pi:G\to\mathcal A$, $(x,h)\mapsto c_xh$, be the
		conjugation map. Since $c$ is $G$-stable, defect-group conjugacy
		gives $G=XM$. Hence $M\cap X=N$, $M/N\cong\bar{G}$ via
		$mN\mapsto mX$, and $\pi(M)=\mathcal A_D$.
		Also $\bar{G}\cong H$ has order prime to $p$.
		The kernel $\mathcal K:=\ker\pi=C_G(X)$ lies in $M$;
		projection onto $H$ gives
		\begin{equation}\label{eq:local-kernel-extension}
			1\longrightarrow Z(X)\longrightarrow\mathcal K
			\longrightarrow H\cap I\longrightarrow1.
		\end{equation}
		
		Choose a maximal $c$-Brauer pair $(D,e_D)$ and put
		\[
		\widetilde E=N_G(D,e_D)/C_X(D),\qquad
		\rho:\widetilde E\to\bar{G},\quad gC_X(D)\mapsto gX
		\quad(g\in N_G(D,e_D)).
		\]
		The maximal-pair Frattini argument makes $\rho$ surjective.
		Since $D\leq C_X(D)$, its kernel is the inertial quotient
		$N_X(D,e_D)/C_X(D)$. Thus $\widetilde E$ is a $p'$-group.
		It acts on $D$ by $a(d)=gdg^{-1}$ for
		$a=gC_X(D)\in\widetilde E$ and $d\in D$; this is well defined
		because $C_X(D)$ centralizes $D$.
		
		The local Clifford class in $H^2(\widetilde E,k^\times)$ has a normalized
		representative with values in roots of unity of order dividing $|\widetilde E|$
		\cite[Proposition~1.2.9(ii)]{Linckelmann2018a}. Let $\kappa$ be its multiplicative lift to
		$\mu_{p'}(\cO)$ as in \cite[p.~6 of the arXiv version]{CoconetMarcusTodea2020}.
		Inflate it to $D\rtimes\widetilde E$ by
		$\kappa((d,a),(d',b))=\kappa(a,b)$ for
		$d,d'\in D$ and $a,b\in\widetilde E$, and define
		\[
		\Lambda:=\cO^\kappa(D\rtimes\widetilde E)
		=\bigoplus_{a\in\widetilde E}\cO D\,\widehat a.
		\]
		Here $\widehat a$ is the standard unit corresponding to $(1,a)$,
		$\widehat1=1$, and we identify $d\in D$ with $d\widehat1$.
		The multiplication and $\bar{G}$-grading are
		\[
		(d\widehat a)(d'\widehat b)
		=\kappa(a,b)d\,a(d')\widehat{ab}
		\qquad(d,d'\in D,\ a,b\in\widetilde E),
		\]
		\[
		\Lambda_{\bar g}
		=\bigoplus_{\substack{a\in\widetilde E\\\rho(a)=\bar g}}
		\cO D\,\widehat a
		\qquad(\bar g\in\bar{G}).
		\]
		Write $\Lambda_1$ for its identity component.
		
		Since $c$ is inertial and $G/X$ has order prime to $p$, choose the graded
		$(\cO Gc,\Lambda)$-Morita bimodule from \cite[Theorems~1.1(c) and~7.8(2)]{CoconetMarcusTodea2020}
		(see also \cite[Proposition~3.3]{Zhou2016}) and the graded $(\Lambda,\cO Mc_0)$-Morita bimodule obtained from the
		local equivalence and Fong--Reynolds correspondence in
		\cite[Corollary~7.7 and Proposition~4.3]{CoconetMarcusTodea2020}. Use the $\bar G$-grading
		induced by $\rho$. In the construction of the first equivalence
		\cite[proof of Theorem~7.8]{CoconetMarcusTodea2020}, the source idempotents $i\in A^D$,
		$i_0\in A_0^D$ and the $N_G(D,e_D)$-invariant endopermutation
		$\cO D$-lattice $V_s$ of rank prime to $p$ give $D$-interior isomorphisms
		\[
		i_0A_0i_0\cong\Lambda_1,\qquad
		\Phi:iAi\xrightarrow{\sim}\End_\cO(V_s)\otimes_\cO\Lambda_1.
		\]
		Write
		$r=\operatorname{rank}_\cO(V_s)$ and let
		$\rho_{V_s}:D\to\GL_\cO(V_s)$ be its representation.
		The structural maps send $di_0$ to $d$ and satisfy
		\begin{equation}\label{eq:source-interior-action}
			\Phi(di)=\rho_{V_s}(d)\otimes d\qquad(d\in D).
		\end{equation}
		At these corners the identity components of the two bimodules are
		$V_s\otimes_\cO\Lambda_1$ and the regular $\Lambda_1$-bimodule, respectively.
		Here the local defect group is $D$ itself, so the local corner over
		$\cO D$ is $\cO D$ \cite[Corollary~7.7]{CoconetMarcusTodea2020}.
		
		We may choose $V_s$ to have determinant $1$
		\cite[p.~446]{LassueurThevenaz2019}. Explicitly, put
		$\delta=\det\rho_{V_s}$ and choose an integer $a$ with
		$ar\equiv1\pmod{\exp(D)}$. Invariance of $V_s$ makes
		$\xi=\delta^{-a}$ invariant under $\widetilde E$. Hence
		\[
		\psi_\xi:\Lambda\longrightarrow\Lambda,\qquad
		d\widehat v\longmapsto\xi(d)d\widehat v
		\quad(d\in D,\ v\in\widetilde E)
		\]
		is a grading-preserving algebra automorphism. Twist the right action on
		the first graded bimodule by $\psi_\xi^{-1}$, so that $w\cdot b=w\psi_\xi^{-1}(b)$ for $b\in\Lambda$.
		On the identity component the map $v\otimes b\mapsto v\otimes\psi_\xi(b)$ identifies
		the twisted right action with the regular one. Under this identification the left action
		is given by $(\id\otimes\psi_\xi|_{\Lambda_1})\Phi$ and sends $di$ to
		$\xi(d)\rho_{V_s}(d)\otimes d$. Thus \eqref{eq:source-interior-action} holds
		with $V_s$ replaced by $\xi\otimes V_s$, whose determinant is $\delta^{1-ar}=1$.
		Retain $V_s,\rho_{V_s},\Phi$ for these normalized data.
		
		For $z\in O_p(Z(X))$, centrality of $zi$ implies that
		$\rho_{V_s}(z)\otimes z$ commutes with every
		$f\otimes1$, $f\in\End_\cO(V_s)$. Thus $\rho_{V_s}(z)$ is
		scalar, say $\rho_{V_s}(z)=\lambda\id_{V_s}$ with
		$\lambda\in\cO^\times$. We have
		$\lambda^{|z|}=1$ and $\lambda^r=\det(\rho_{V_s}(z))=1$.
		Since $(|z|,r)=1$, it follows that
		$\rho_{V_s}(z)=\id_{V_s}$. Consequently
		$\Phi(zi)=\id_{V_s}\otimes z$, and $z$ is central in $\Lambda_1$.
		
		Using the two identity-component isomorphisms, form the
		$(A,A_0)$-Morita bimodule
		\[
		U=Ai\otimes_{iAi}(V_s\otimes_\cO\Lambda_1)
		\otimes_{\Lambda_1}i_0A_0.
		\]
		The source idempotents are full, and the middle factor is the
		usual matrix Morita bimodule. The preceding calculation shows
		that $zu=uz$ for $z\in O_p(Z(X))$. Brauer correspondence
		matches the central characters on $O_{p'}(Z(X))$, so the same
		identity holds there. Thus $U$ is balanced over $\cO Z(X)$.
		
		Tensor the twisted first graded bimodule with the local one over $\Lambda$.
		This gives a $\bar G$-graded $(\cO Gc,\cO Mc_0)$-Morita bimodule $W$.
		Strong grading identifies its identity component with the tensor product
		of the identity components over $\Lambda_1$, hence with $U$.
		Use this identification in what follows.
		
		For $g\in M$, define $t_g(u)=gug^{-1}$ for $u\in U$,
		using the two actions on $W$. Conjugation preserves the
		identity degree, and the bimodule identities give an
		$\cO$-linear action. For $g,h\in M$, $a\in A$, $b\in A_0$
		and $u\in U$,
		\[
		t_gt_h=t_{gh},\qquad
		t_g(aub)=(gag^{-1})t_g(u)(gbg^{-1}).
		\]
		Central balance gives $t_z=\id_U$ for $z\in Z(X)$.
		Since $\mathcal K$ centralizes $X$, its image lies in
		$\Aut_{A\text{-}A_0}(U)\cong Z(A)^\times$.
		Since $t(Z(X))=1$, this image factors through
		$\mathcal K/Z(X)\cong H\cap I$. It is therefore finite of order
		prime to $p$, by \eqref{eq:local-kernel-extension}.
		The ring $Z(A)$ is complete local with residue field $k$.
		Uniqueness in Hensel's lemma identifies each of its
		roots of unity of order prime to $p$ with the scalar lift
		of its residue in $\cO$. We therefore identify
		\[
		\mu:=t(\mathcal K)\leq\mu_{p'}(\cO).
		\]
		
		Since $\pi(M)=\mathcal A_D$ and $\pi(N)=I_D$,
		we have $M/(N\mathcal K)\cong\Omega_D$.
		Choose representatives $m_\omega\in M$ for $\omega\in\Omega_D$,
		with $m_1=1$. For $\omega,\eta\in\Omega_D$, choose
		$n_{\omega,\eta}\in N$ and $k_{\omega,\eta}\in\mathcal K$ such that
		\[
		m_\omega m_\eta=n_{\omega,\eta}k_{\omega,\eta}m_{\omega\eta},
		\qquad t_{k_{\omega,\eta}}=\beta(\omega,\eta)\id_U.
		\]
		This defines $\beta(\omega,\eta)\in\mu$ independently of the
		factorization, since $N\cap\mathcal K=Z(X)$ and $t(Z(X))=1$.
		
		For $\alpha\in\mathcal A_D$, put $\omega=\alpha I_D$,
		write $\alpha=c_n\pi(m_\omega)$ with $n\in N$, and define
		\[
		T_\alpha=t_n t_{m_\omega}.
		\]
		The choice of $n$ is unique modulo $Z(X)$, so $T_\alpha$
		is well defined and semilinear for $\alpha$.
		For $\alpha=c_n\pi(m_\omega)$ and
		$\alpha'=c_{n'}\pi(m_\eta)$, with $n,n'\in N$, the group law gives
		\[
		\begin{aligned}
			T_\alpha T_{\alpha'}
			&=\beta(\omega,\eta)\,
			t_{n\,m_\omega n'm_\omega^{-1}n_{\omega,\eta}}t_{m_{\omega\eta}}\\
			&=\beta(\omega,\eta)T_{\alpha\alpha'}.
		\end{aligned}
		\]
		As $m_1=1$, we have $T_1=\id_U$ and
		$T_{c_n}(u)=nun^{-1}$ for $n\in N$ and $u\in U$.
		The multiplication law and associativity show that
		$\beta\in Z^2(\Omega_D,\mu)$ is normalized.
		Together with central balance, these identities show that
		$(U,\mu,\beta,\{T_\alpha\})$ is a scalar-equivariant Morita datum
		for $(X,c,D)$ in the sense of Definition~\ref{def:datum}.
	\end{proof}
	
	\subsection{Removal of the block kernel}
	
	\begin{proof}[Proof of Proposition~\ref{prop:local}]
		The group $\ker(c)$ is a normal $p'$-subgroup of $X$
		\cite[Theorem~6.10]{Navarro1998}.
		If $\ker(c)=X$, then $D=1$, $N=X$, and $c_0=c$;
		the identity bimodule $\cO Xc$ with its natural
		$\Aut(X)_c$-action gives the required maps with multiplier $1$.
		Otherwise $\ker(c)\leq Z(X)$. For $z\in\ker(c)$, we have $zc=c$
		by definition and $zc_0=c_0$ by Brauer correspondence.
		Put $Z=\ker(c)$ and write $\bar X=X/Z$ and $\bar x=xZ$
		for $x\in X$. Set $\bar D=DZ/Z$ and
		$\bar N=N/Z=N_{\bar X}(\bar D)$,
		and let $\bar c,\bar c_0$ be the quotient blocks.
		The normalizer equality holds because the preimage of $\bar D$
		is $DZ=D\times Z$, whose unique Sylow $p$-subgroup is $D$.
		The quotient map induces algebra isomorphisms
		\[
		\varphi:\cO Xc\xrightarrow{\sim}\cO\bar X\bar c,\qquad
		\varphi_0:\cO Nc_0\xrightarrow{\sim}\cO\bar N\bar c_0,
		\]
		given by $\varphi(xc)=\bar x\bar c$ and
		$\varphi_0(nc_0)=\bar n\bar c_0$ for $x\in X$ and $n\in N$.
		These are compatible with Brauer correspondence and inflation
		of characters \cite[Theorem~9.9(c)]{Navarro1998}.
		
		Every $\alpha\in\Aut(X)_c$ stabilizes $Z$ and induces
		$\bar\alpha\in\Aut(\bar X)_{\bar c}$ by
		$\bar\alpha(\bar x)=\overline{\alpha(x)}$ for $x\in X$.
		For $\alpha\in\mathcal A_D=\Aut(X)_{D,c}$, the induced
		automorphism $\bar\alpha$ stabilizes $\bar D$, and
		\[
		\varphi\circ\alpha=\bar\alpha\circ\varphi,\qquad
		\varphi_0\circ\alpha=\bar\alpha\circ\varphi_0.
		\]
		Put $I_D=\{c_n:n\in N\}$, $\Omega_D=\mathcal A_D/I_D$, and
		\[
		\bar I_{\bar D}=\{c_{\bar n}:\bar n\in\bar N\},\qquad
		\bar\Omega_{\bar D}=
		\Aut(\bar X)_{\bar D,\bar c}/\bar I_{\bar D}.
		\]
		Since $\overline{c_n}=c_{\bar n}$ for $n\in N$, the quotient action
		induces a homomorphism
		\[
		\eta:\Omega_D\longrightarrow\bar\Omega_{\bar D},\qquad
		\eta(\alpha I_D)=\bar\alpha\bar I_{\bar D}
		\quad(\alpha\in\mathcal A_D).
		\]
		
		The group $\bar X$ is quasisimple, $\bar c$ is faithful, and
		equivariant inflation makes $\Aut(\bar X)_{\bar c}$ transitive
		on $\IBr(\bar c)$.
		Proposition~\ref{prop:classification} and
		Lemmas~\ref{lem:nilpotent-covering} and~\ref{lem:outer-pprime}
		therefore give a scalar-equivariant Morita datum
		$(\bar U,\mu,\bar\beta,\{\bar T_\gamma\})$ for
		$(\bar X,\bar c,\bar D)$ in the sense of Definition~\ref{def:datum}.
		In particular, $\bar D$ is abelian. Since $D\cap Z=1$,
		the group $D\cong\bar D$ is abelian as well.
		Give $U=\bar U$ the actions $aub=\varphi(a)u\varphi_0(b)$
		for $a\in\cO Xc$ and $b\in\cO Nc_0$, and set
		\[
		T_\alpha=\bar T_{\bar\alpha}\quad(\alpha\in\mathcal A_D),
		\qquad \beta=\eta^*\bar\beta.
		\]
		The intertwining identities give semilinearity, while
		$\overline{c_n}=c_{\bar n}$ and $Z(X)/Z\leq Z(\bar X)$
		ensure the prescribed inner action and central balance.
		The identities in \eqref{eq:scalar-datum} follow by pullback.
		Thus $(U,\mu,\beta,\{T_\alpha\})$ is a scalar-equivariant Morita datum
		for $(X,c,D)$ in the sense of Definition~\ref{def:datum}.
	\end{proof}

	\section{Proof of Theorem~\ref*{thm:solvable}}\label{sec:reduction}
	
	\begin{proof}[Proof of Theorem~\ref{thm:solvable}]
		For a finite group $J$, put
		\[
		Z_0(J)=O_{p'}(Z(J)),\qquad m(J)=|J:Z_0(J)|.
		\]
		Suppose that the theorem is false, and choose a counterexample
		$(G,b)$ with $m(G)$ minimal. All Morita equivalences below are
		over the fixed coefficient ring $\cO$.
		
		Let $H\unlhd G$, let $d$ be a block of $\cO H$ covered by $b$,
		and put $I=I_G(d)$. Fong--Reynolds reduction
		\cite[Theorem~6.8.3]{Linckelmann2018b} gives a block $b_I$ of $\cO I$
		Morita equivalent to $b$. In particular, $l(b_I)=1$.
		Since $Z_0(G)\leq O_{p'}(Z(I))$, if $I<G$ then
		\[
		m(I)\leq |I:Z_0(G)|<m(G).
		\]
		Minimality would make $\cO Ib_I$, and hence $\cO Gb$,
		Morita equivalent to a block algebra of a finite solvable
		group, a contradiction. Thus $I=G$, and every block of a normal
		subgroup covered by $b$ is $G$-stable. Since the covered blocks of a fixed normal
		subgroup form one $G$-orbit \cite[Corollary~9.3]{Navarro1998}, there is a unique such block.
		
		Suppose that $H\unlhd G$ and that $b$ covers a nilpotent block
		$d_H$ of $\cO H$. Put $H_1=HZ(G)$, and choose a block $d$ of
		$\cO H_1$ covered by $b$ and covering $d_H$. By the preceding
		argument, $d$ is $G$-stable. It is nilpotent because adjoining a central
		subgroup gives a central product, and a block of that product
		covering the nilpotent block $d_H$ is nilpotent \cite[Lemma~7.5(iii)]{Sambale2014}.
		Choose a defect group $P_1$ of $b$ such that
		$Q=P_1\cap H_1$ is a defect group of $d$
		\cite[Theorem~9.26]{Navarro1998}.
		The K\"ulshammer--Puig extension theorem \cite[Theorems~1.8 and~1.12]{KulshammerPuig1990}, in the
		form recalled in \cite[Lemma~5.2]{JiangZhangZhou2026}, gives a finite group
		$G'$, a central $p'$-subgroup $Z'\leq G'$, and a block $b'$
		of $\cO G'$ Morita equivalent to $b$, with
		\[
		|G':Z'|=|G:H_1|\,|Q|.
		\]
		As $QZ_0(G)\leq H_1$ and $Q\cap Z_0(G)=1$, we obtain
		\[
		m(G')\leq |G':Z'|=|G:H_1|\,|Q|\leq m(G).
		\]
		The pair $(G',b')$ is also a counterexample: $l(b')=1$,
		and a Morita equivalence of $\cO G'b'$ with a block algebra
		of a finite solvable group would give one for $\cO Gb$.
		Minimality therefore forces equality throughout.
		In particular, $|H_1|=|Q||Z_0(G)|$, so
		$H_1=Q\times Z_0(G)$. Its unique Sylow $p$-subgroup $Q$
		is characteristic in $H_1$ and hence normal in $G$.
		Consequently,
		\begin{equation}\label{eq:nilpotent-normal-subgroup}
			H\leq O_p(G)Z_0(G)
		\end{equation}
		whenever $H\unlhd G$ and $b$ covers a nilpotent block of
		$\cO H$.
		
		Applying \eqref{eq:nilpotent-normal-subgroup} to $O_{p'}(G)$,
		whose covered block has defect zero, gives
		\[
		O_{p'}(G)=Z_0(G).
		\]
		Recall that the Fitting subgroup $F(G)$ is the largest
		nilpotent normal subgroup, that a component is a subnormal
		quasisimple subgroup, and that the layer $E(G)$ is the
		product of the components. The generalized Fitting subgroup
		$F^*(G)=F(G)E(G)$ is self-centralizing \cite[Theorem~9.8]{Isaacs2008}.
		If $E(G)=1$, put $Q=O_p(G)$ and $Z=Z_0(G)$. Then
		\[
		F^*(G)=Q\times Z,\qquad
		C_G(Q)\leq QZ,\qquad C_G(Q)=Z(Q)\times Z.
		\]
		Let $\lambda\in\Irr(Z)$ be the unique irreducible character
		satisfying
		\[
		\chi|_Z=\chi(1)\lambda\qquad(\chi\in\Irr(b)).
		\]
		Since $Z$ is a $p'$-group, we also regard $\lambda$ as a
		Brauer character. For its primitive idempotent $e_\lambda\in kZ$,
		we have $be_\lambda=b$ and
		\[
		kC_G(Q)e_\lambda\cong kZ(Q),
		\]
		which is local.
		Decompose $e_\lambda$ as a sum of block idempotents of $kG$.
		Every summand $d$ satisfies $de_\lambda=d$. Since $Q$ is a
		normal $p$-subgroup, it lies in a defect group of every block;
		hence $\Br_Q(d)\ne0$ \cite[Theorem~9.26]{Navarro1998}. These Brauer images
		are orthogonal idempotents in $kC_G(Q)e_\lambda$, with sum
		$\Br_Q(e_\lambda)=e_\lambda$. A local algebra has only the
		idempotents zero and its identity, so there is just one summand.
		Thus $e_\lambda$ is a block idempotent of $kG$ and $b=e_\lambda$.
		Writing $\IBr(G\mid\lambda)$ for the irreducible Brauer
		characters whose restrictions to $Z$ contain $\lambda$,
		we obtain
		\[
		\IBr(G\mid\lambda)=\IBr(b)=\{\varphi\},\qquad
		\lambda^G=t\varphi
		\]
		for some positive integer $t$. By \cite[Theorem~1]{Sambale2024},
		$G/Z$ is solvable. Since $Z$ is abelian, $G$ is solvable,
		contrary to the choice of $(G,b)$. Therefore $E(G)\ne1$.
		
		Choose one $G$-orbit of components $X_1,\ldots,X_r$, and put
		$L=X_1\cdots X_r$. Choose a block $e$ of $\cO E(G)$ covered
		by $b$, and let $c$ be a block of $\cO L$ covered by $e$.
		Then $b$ covers $c$, so $c$ is $G$-stable.
		Choose a defect group $P$ of $b$ such that
		$D=P\cap L$ is a defect group of $c$
		\cite[Theorem~9.26]{Navarro1998}.
		For each $i$, let $f_i$ be the block of $\cO X_i$ covered by
		$c$, and put $D_i=D\cap X_i$. The block theory of central
		products gives
		\[
		D=D_1\cdots D_r,
		\]
		where $D_i$ is a defect group of $f_i$ \cite[Lemma~7.5(i)]{Sambale2014}.
		As $l(b)=1$, the set $\IBr(b)$ is a single
		$\Aut(G)_b$-orbit. Each $f_i$ is covered by the block $e$ of
		$E(G)$, itself covered by $b$. Thus \cite[Lemma~2.5]{JiangZhangZhou2026}
		applies, even though the components need not be normal in $G$,
		and shows that $\Aut(X_i)_{f_i}$ is transitive on $\IBr(f_i)$.
		Proposition~\ref{prop:local} therefore gives that every $D_i$
		is abelian. Since distinct components commute, $D$ is abelian.
		
		If $D\leq Z(L)$, then $c$ is nilpotent, so
		\eqref{eq:nilpotent-normal-subgroup} gives
		$L\leq O_p(G)Z_0(G)$. This is impossible because $L$ is a
		nontrivial perfect group. Thus $D\not\leq Z(L)$.
		As $L\unlhd G$, we have $D=P\cap L\unlhd P$, whence
		$P\leq N_G(D)$. If $D\unlhd G$, each $D_i$ is a normal
		$p$-subgroup of the quasisimple group $X_i$, and hence
		$D_i\leq Z(X_i)$. This would give $D\leq Z(L)$, a
		contradiction. We have proved
		\[
		P\leq M:=N_G(D)<G.
		\]
		Put $N=N_L(D)$, and let $c_0$ be the Brauer correspondent
		of $c$ in $N$. Since $c$ is $G$-stable, defect-group
		conjugacy gives $G=LM$ and $M/N\cong G/L$.
		
		By Proposition~\ref{prop:local}, $(X_1,f_1,D_1)$ admits a
		scalar-equivariant Morita datum with scalar group $\mu$.
		Here $L$ is the product of the $G$-conjugates of $X_1$,
		the block $c$ is $G$-stable and covered by $b$, and
		$D_1=D\cap X_1$. Thus all hypotheses of
		Proposition~\ref{prop:quotient-induction} hold. Applied to $b$,
		$X_1$, and $c$, it gives a central extension
		\[
		1\longrightarrow\mu\longrightarrow\widehat M
		\xrightarrow{\rho}M\longrightarrow1
		\]
		and a block $b'$ of $\cO\widehat M$ such that
		\[
		\cO Gb\sim_{\Mor}\cO\widehat M b',\qquad l(b')=l(b)=1.
		\]
		
		The same proposition gives
		$\rho^{-1}(Z_0(G))=\mu\times Z_0(G)\leq Z(\widehat M)$.
		Thus $\mu\times Z_0(G)$ is a central $p'$-subgroup of $\widehat M$.
		As $M<G$, we obtain
		\begin{align}
			m(\widehat M)
			&=|\widehat M:O_{p'}(Z(\widehat M))|\nonumber\\
			&\leq|\widehat M:\mu\times Z_0(G)|
			=|M:Z_0(G)|<|G:Z_0(G)|=m(G).\label{eq:strict-central-index}
		\end{align}
		By minimality, $\cO\widehat M b'$ is Morita equivalent over $\cO$
		to a block algebra of a finite solvable group. Transitivity of Morita
		equivalence shows that $\cO Gb$ is Morita equivalent over $\cO$ to a
		block algebra of a finite solvable group, a contradiction.
	\end{proof}

	\section{A derivation criterion for twisted group algebras}\label{sec:derivations}
	
	Throughout this section $k$ is algebraically closed of characteristic $p$.
	Let $A=k^\alpha L$ be a twisted group algebra of a finite group $L$.
	Put $N=O_p(L)$, $H=L/N$, and $\bar g=gN$ for $g\in L$.
	Since $N$ is a $p$-group, $\operatorname{Hom}(N,k^\times)=1$ and
	$H^2(N,k^\times)=1$ \cite[Proposition~1.2.9(iii)]{Linckelmann2018a}. Thus inflation
	$H^2(H,k^\times)\to H^2(L,k^\times)$ is an isomorphism
	\cite[Chapter~VII, \S6, Proposition~5]{Serre1979}.
	Put $Z=\{z\in k^\times:z^{|H|}=1\}$. By
	\cite[Proposition~1.2.9(ii)]{Linckelmann2018a}, choose a normalized $2$-cocycle
	$\beta:H\times H\to Z$ whose inflation represents $[\alpha]$.
	We take this inflation as the representative $\alpha$ and rescale the
	basis $(u_g)_{g\in L}$ of $A$ accordingly. Then $u_1=1$ and
	\[
	u_gu_h=\alpha(g,h)u_{gh},\qquad
	\alpha(g,h)=\beta(\bar g,\bar h)\quad(g,h\in L).
	\]
	Write $(w_x)_{x\in H}$ for the basis of $k^\beta H$ with
	$w_xw_y=\beta(x,y)w_{xy}$. Then $u_g\mapsto w_{\bar g}$ defines a
	surjective algebra homomorphism $A\to k^\beta H$.
	
	A derivation of a $k$-algebra is a $k$-linear map $\partial$ from the
	algebra to itself satisfying $\partial(xy)=x\partial(y)+\partial(x)y$.
	
	\begin{proposition}[Derivation criterion]\label{prop:derivation}
		Let $O_{p'}(L)=1$ and suppose that $A=k^\alpha L$ has just one simple
		module up to isomorphism. Let $B$ be its basic algebra and $J=J(B)$.
		The following conditions are equivalent.
		\begin{enumerate}
			\item $L$ is a $p$-group.
			\item For each $x\in J\setminus J^2$ there is a derivation $\partial$ of $B$
			such that $\partial(x)\notin J$.
		\end{enumerate}
	\end{proposition}
	
	Here $B=eAe$ for a primitive idempotent $e$ of $A$. It is local,
	with residue field $k$. Fix a representative $V$ of the unique isomorphism
	class of simple $A$-modules. Set $\mathcal{E}=\End_k(V)$, and let $\rho:A\to \mathcal{E}$
	be the representation afforded by $V$. The nonzero quotient $k^\beta H$
	also has a unique simple module up to isomorphism, whose inflation is
	isomorphic to $V$. Thus $u_g$ acts on $V$ as $w_{\bar g}$, and
	$\rho(u_n)=\mathrm{id}_V$ for $n\in N$ since $w_1=1$.
	We first interpret condition~(ii)
	on $V$ and then use group transfer to test it on $N$.
	
	\subsection{Derivations and the simple module}
	
	Let $\varepsilon:B\to B/J=k$ be the quotient map. For a derivation $\partial$
	of $B$ and $x,y\in J$, the Leibniz rule gives
	\[
	(\varepsilon\circ\partial)(xy)
	=\varepsilon(\partial(x))\varepsilon(y)+\varepsilon(x)\varepsilon(\partial(y))=0,
	\]
	since $\varepsilon(x)=\varepsilon(y)=0$. Such products span $J^2$, so the
	restriction of $\varepsilon\circ\partial$ to $J$ induces a $k$-linear map $J/J^2\to k$.
	These functionals form a subspace of $(J/J^2)^*$. Condition~(ii) of
	Proposition~\ref{prop:derivation} says that their common annihilator in $J/J^2$ is zero.
	By finite dimensionality, this is equivalent to their spanning the whole dual.
	
	We call a $k$-linear map $\delta:A\to \mathcal{E}$ a \emph{$\rho$-derivation} if
	\[
	\delta(xy)=\rho(x)\delta(y)+\delta(x)\rho(y)\qquad(x,y\in A).
	\]
	The products on the right are compositions in $\mathcal{E}=\End_k(V)$.
	
	\begin{lemma}\label{lem:lift-derivation}
		Let $A=k^\alpha L$, let $V$ be its unique simple module up to isomorphism,
		and let $\rho:A\to \mathcal{E}=\End_k(V)$ be the corresponding representation.
		Condition~(ii) of Proposition~\ref{prop:derivation} for the basic algebra $B$
		holds if and only if every $\rho$-derivation $\delta:A\to \mathcal{E}$ has the form
		$\rho\circ\partial$ for a derivation $\partial:A\to A$.
	\end{lemma}
	
	\begin{proof}
		Put $d=\dim_k V$. Since $V$ is the only simple $A$-module up to
		isomorphism, we may write $A$ as the algebra of $d\times d$ matrices
		over $B$ \cite[Proposition~4.7.17(i), (ii), (iv)]{Linckelmann2018a}.
		Recall that $\varepsilon:B\to B/J\simeq k$ is the quotient map.
		For each $b\in B$, $\varepsilon(b)$ is the unique scalar $a\in k$
		such that $b-a1_B\in J$. In a suitable basis of $V$, the matrix
		$(b_{ij})$ therefore acts as $(\varepsilon(b_{ij}))$.
		Write $I$ for the identity matrix and regard
		matrices over $k$ as matrices over $B$ using $k1_B\subseteq B$.
		For $b\in B$ and $S=(s_{ij})$ over $k$, write $bS=(bs_{ij})$.
		
		First suppose that condition~(ii) of Proposition~\ref{prop:derivation} holds,
		and let $\delta:A\to \mathcal{E}$ be a $\rho$-derivation. For $b\in B$ and any matrix $S$ over
		$k$, we have $(bI)S=S(bI)$. Applying $\delta$ to both sides and using
		$\rho(bI)=\varepsilon(b)I$ gives $\delta(bI)S=S\delta(bI)$.
		Only scalar matrices commute with all matrices over $k$, so
		$\delta(bI)=f(b)I$ for a $k$-linear map $f:B\to k$. The product rule gives
		\[
		f(bc)=\varepsilon(b)f(c)+f(b)\varepsilon(c).
		\]
		Thus $f(1_B)=0$ and $f(J^2)=0$. By the paragraph before the lemma,
		condition~(ii) gives a derivation $\partial_B$ of $B$ with
		$f=\varepsilon\circ\partial_B$, since $B=k1_B\oplus J$.
		
		We now construct the required derivation of $A$. Set
		\[
		\partial(bS)=\partial_B(b)S+b\,\delta(S).
		\]
		To define $\partial$ on any matrix, write it as a sum of matrices with
		just one nonzero entry and apply this formula to each summand.
		For $b,c\in B$ and matrices $S,T$ over $k$, we have $(bS)(cT)=bc(ST)$.
		Since $\rho$ fixes matrices over $k$, we have
		$\delta(ST)=\delta(S)T+S\delta(T)$. Hence
		\[
		\begin{aligned}
			\partial((bS)(cT))
			&= \partial_B(bc)ST+bc\,\delta(ST)\\
			&= \bigl(\partial_B(b)c+b\partial_B(c)\bigr)ST
			+bc\bigl(\delta(S)T+S\delta(T)\bigr)\\
			&= \partial(bS)(cT)+(bS)\partial(cT).
		\end{aligned}
		\]
		Matrices of the form $bS$ span $A$ over $k$, so bilinearity gives
		the same identity for all pairs of matrices in $A$. Thus $\partial$
		is a derivation of $A$.
		Moreover,
		\[
		\rho(\partial(bS))=f(b)S+\varepsilon(b)\delta(S)=\delta((bI)S).
		\]
		Thus $\rho\circ\partial=\delta$.
		
		Conversely, suppose every $\delta$ in the statement has such a lift.
		Given $\lambda\in\Hom_k(J/J^2,k)$, put $f(a1_B+x)=\lambda(x+J^2)$
		for $a\in k$ and $x\in J$. To compute $f(bc)$, write
		$b=r1_B+x$ and $c=s1_B+y$, with $r,s\in k$ and $x,y\in J$. Then
		\[
		bc=rs1_B+ry+sx+xy.
		\]
		The definition of $f$ discards the scalar term $rs1_B$. Also,
		$xy\in J^2$, so $xy+J^2=0$ in $J/J^2$. By linearity of $\lambda$,
		\[
		\begin{aligned}
			f(bc)&=\lambda(ry+sx+xy+J^2)\\
			&=r\lambda(y+J^2)+s\lambda(x+J^2)\\
			&=\varepsilon(b)f(c)+f(b)\varepsilon(c),
		\end{aligned}
		\]
		where the last equality uses $\varepsilon(b)=r$ and $\varepsilon(c)=s$.
		For $X=(x_{ij})\in A$, define $\delta(X)=(f(x_{ij}))$.
		If $Y=(y_{ij})\in A$, then $(XY)_{ij}=\sum_{\ell=1}^d x_{i\ell}y_{\ell j}$.
		By linearity of $f$ and the formula for $f(bc)$ above,
		\[
		\begin{aligned}
			(\delta(XY))_{ij}
			&=\sum_{\ell=1}^d f(x_{i\ell}y_{\ell j})\\
			&=\sum_{\ell=1}^d \varepsilon(x_{i\ell})f(y_{\ell j})
			+\sum_{\ell=1}^d f(x_{i\ell})\varepsilon(y_{\ell j})\\
			&=(\rho(X)\delta(Y))_{ij}+(\delta(X)\rho(Y))_{ij}.
		\end{aligned}
		\]
		Since this holds for all $i,j$, we have
		$\delta(XY)=\rho(X)\delta(Y)+\delta(X)\rho(Y)$. Choose a derivation $\partial$
		of $A$ with $\rho\circ\partial=\delta$. For $b\in B$, define
		$\partial_B(b)=(\partial(bI))_{11}$, the entry in the first row and first
		column of the matrix $\partial(bI)$. Comparing the $(1,1)$ entries in
		$\partial((bI)(cI))=(bI)\partial(cI)+\partial(bI)(cI)$ gives
		$\partial_B(bc)=b\partial_B(c)+\partial_B(b)c$. Thus $\partial_B$ is
		a derivation of $B$, and $\varepsilon(\partial_B(x))=\lambda(x+J^2)$
		for every $x\in J$. This proves condition~(ii) of Proposition~\ref{prop:derivation}.
	\end{proof}
	
	\begin{lemma}\label{lem:one-simple-reduction}
		Suppose that $O_{p'}(L)=1$ and that $k^\alpha L$ has only one simple
		module up to isomorphism. Then $L$ is solvable, $C_L(N)\leq N$, and
		$O_p(H)=1$. Conjugation induces a faithful action of $H$ on
		$N/\Phi(N)$, where $\Phi(N)=[N,N]N^p$ is the Frattini subgroup of $N$ \cite[Proposition~2.25]{Craven2008}.
	\end{lemma}
	
	\begin{proof}
		Since $\alpha(g,h)=\beta(\bar g,\bar h)\in Z$, the set
		$\widehat L=\{zu_g:z\in Z,\ g\in L\}$ is a finite group under
		multiplication in $A$. The map $zu_g\mapsto g$ has central kernel $Z$.
		Let $\lambda:Z\hookrightarrow k^\times$ be the inclusion. The section
		$g\mapsto u_g$ has factor set $\alpha$, since
		$u_gu_hu_{gh}^{-1}=\alpha(g,h)\in Z$. Thus
		\cite[Proposition~1.2.17]{Linckelmann2018a} gives a surjective algebra homomorphism
		\[
		\pi:k\widehat L\longrightarrow A,\qquad zu_g\longmapsto\lambda(z)u_g.
		\]
		Let $k_\lambda$ denote the one-dimensional left $kZ$-module with
		$k$-basis $\{v_\lambda\}$ on which each $z\in Z$ acts by
		$z\cdot v_\lambda=\lambda(z)v_\lambda$.
		Since $\pi(z)=\lambda(z)1_A$, the map $x\otimes v_\lambda\mapsto\pi(x)$
		gives a surjection from $\operatorname{Ind}_Z^{\widehat L}k_\lambda
		=k\widehat L\otimes_{kZ}k_\lambda$ onto the regular left $A$-module,
		viewed as a $k\widehat L$-module via $\pi$. It is an isomorphism
		because both modules have dimension $|\widehat L:Z|=|L|$.
		As $\pi$ is surjective, the $k\widehat L$-submodules of $A$ are exactly
		its left ideals. Thus the induced module has the same composition factors
		as the regular $A$-module, all isomorphic to $V$ by hypothesis.
		Hence the Brauer character of
		$k_\lambda$ induces to a positive multiple of the Brauer character of
		$V$. By Sambale's theorem \cite[Theorem~1]{Sambale2024}, $L\cong\widehat L/Z$ is solvable.
		
		Since $O_{p'}(L)=1$, the Fitting subgroup of $L$ is $N$. Thus
		$C_L(N)\leq N$ by the centralizer theorem for the Fitting subgroup of a
		solvable group \cite[Theorem~9.8]{Isaacs2008}.
		The inverse image in $L$ of any normal $p$-subgroup of $H$ is a normal
		$p$-subgroup of $L$, so $O_p(H)=1$.
		
		Let $\theta:L\to\Aut(N/\Phi(N))$ be the homomorphism induced by
		conjugation. Since $[N,N]\leq\Phi(N)$, we have $N\leq\ker\theta$.
		By the Hall--Burnside theorem \cite[Theorem~2.30]{Craven2008},
		$(\ker\theta)/C_L(N)$ is a $p$-group. The subgroup $C_L(N)$ of $N$
		is a $p$-group as well. Hence $\ker\theta$ is a normal $p$-subgroup of $L$.
		Thus $\ker\theta\leq O_p(L)=N$, so $\ker\theta=N$ and $H=L/N$
		acts faithfully on $N/\Phi(N)$.
	\end{proof}
	
	Recall from the paragraph following Proposition~\ref{prop:derivation}
	that $\mathcal{E}=\End_k(V)$ and that $\rho:A\to \mathcal{E}$ is the representation on $V$.
	There we saw that $u_g$ acts on $V$ as $w_{gN}$, so $\rho(u_g)$ depends
	only on $gN$. We now define an action of $H=L/N$ on $\mathcal{E}$ by
	\[
	(gN)\cdot a=\rho(u_g)a\rho(u_g)^{-1}\qquad(g\in L,\ a\in \mathcal{E}).
	\]
	For $g,h\in L$, the operators $\rho(u_g)\rho(u_h)$ and $\rho(u_{gh})$
	differ by the nonzero scalar $\beta(gN,hN)$ and hence induce the same
	conjugation on $\mathcal{E}$. This makes $\mathcal{E}$ a $kH$-module.
	The next lemma constructs a $\rho$-derivation, as defined before
	Lemma~\ref{lem:lift-derivation}, that is nonzero on some $u_n$ with $n\in N$.
	
	\begin{lemma}\label{lem:target-derivation}
		Suppose that $O_{p'}(L)=1$ and that $A=k^\alpha L$ has a unique simple
		module $V$ up to isomorphism. Put $N=O_p(L)$ and $H=L/N$, and let
		$q:L\to H$ be the quotient map. Let $T$ be a nontrivial normal
		$p'$-subgroup of $H$. Choose $T_0\leq q^{-1}(T)$ such that
		$q^{-1}(T)=N\rtimes T_0$, and put $P_0=[N,T_0]$. Set
		\[
		e_T=\frac{1}{|T|}\sum_{t\in T}t\in kH,\qquad \mathcal{E}_1=(1-e_T)\mathcal{E}.
		\]
		There exists a $\rho$-derivation $\delta:A\to \mathcal{E}$ such that
		$\delta(u_n)\in\mathcal{E}_1$ for every $n\in N$ and
		$\delta(u_n)\ne0$ for at least one $n\in P_0$.
	\end{lemma}
	
	\begin{proof}
		The group $q^{-1}(T)=\{g\in L:gN\in T\}$ has $N$ as a normal Sylow
		$p$-subgroup, since $q^{-1}(T)/N\cong T$ has order prime to $p$.
		Schur--Zassenhaus gives the complement $T_0$, and $q$ restricts to an
		isomorphism $T_0\to T$.
		
		Averaging over $T$ gives a $k$-linear projection $a\mapsto e_Ta$
		with image $e_T\mathcal{E}=\mathcal{E}^T$. Since $T\unlhd H$,
		conjugation by each $h\in H$ permutes $T$, so $he_T=e_Th$ in $kH$.
		Thus $e_T(h\cdot a)=h\cdot(e_Ta)$ for all $h\in H$ and
		$a\in\mathcal{E}$, so the projection $\mathcal{E}\to\mathcal{E}^T$
		is a homomorphism of $kH$-modules.
		It follows that $\mathcal{E}=\mathcal{E}^T\oplus\mathcal{E}_1$ as $kH$-modules
		and $(\mathcal{E}_1)^T=0$.
		
		Put $X=k\otimes_{\Fp}(N/\Phi(N))$, with the action induced by
		conjugation, and let $X_1=(1-e_T)X$. By Lemma~\ref{lem:one-simple-reduction},
		$H$ acts faithfully on $N/\Phi(N)$. Since $T\ne1$, some element of $T$
		acts nontrivially on $X$. Hence $e_TX=X^T\ne X$ and $X_1\ne0$.
		Give $X_1\otimes_k V$ the diagonal $k^\beta H$-action
		\cite[Proposition~1.6.10]{Linckelmann2018a} given by
		\[
		w_h(x\otimes v)=(h\cdot x)\otimes(w_hv)
		\qquad(h\in H,\ x\in X_1,\ v\in V).
		\]
		This nonzero finite-dimensional module
		has a simple quotient. As $V$ is the unique simple $k^\beta H$-module
		up to isomorphism, there is a surjective homomorphism of $k^\beta H$-modules
		$\psi:X_1\otimes_k V\to V$. For each $x\in X_1$, define the
		endomorphism $F(x)\in\mathcal{E}=\End_k(V)$ by
		\[
		F(x)(v)=\psi(x\otimes v)\qquad(v\in V).
		\]
		The resulting $k$-linear map $F:X_1\to\mathcal{E}$ is nonzero since $\psi\ne0$.
		Since $\psi$ is a module homomorphism, the diagonal action gives
		$\psi((h\cdot x)\otimes w_hv)=w_h\psi(x\otimes v)$ for every $h\in H$.
		Hence $F(h\cdot x)(w_hv)=w_hF(x)(v)$. The operator $w_h$ on $V$ is
		invertible, so this identity says $F(h\cdot x)=h\cdot F(x)$ for the
		conjugation action on $\mathcal{E}$ defined above. Hence $F(e_Tx)=e_TF(x)$.
		For $x\in X_1$ we have $e_Tx=0$, and therefore $F(x)\in\mathcal{E}_1$.
		
		Define $f:N\to(\mathcal{E}_1,+)$ by
		\[
		f(n)=F\bigl((1-e_T)(1\otimes(n\Phi(N)))\bigr).
		\]
		Since $\Phi(N)=[N,N]N^p$, the quotient $N/\Phi(N)$ is abelian and
		every element has $p$th power equal to the identity. We write this
		$\Fp$-vector space additively, so for $n,m\in N$ the product of
		their cosets is written as $(nm)\Phi(N)=n\Phi(N)+m\Phi(N)$.
		Consequently,
		\[
		1\otimes((nm)\Phi(N))=1\otimes(n\Phi(N))+1\otimes(m\Phi(N)).
		\]
		Applying the $k$-linear maps $1-e_T$ and $F$ gives
		$f(nm)=f(n)+f(m)$. Thus $f:N\to(\mathcal{E}_1,+)$ is a group
		homomorphism. Equivariance gives
		$f(gng^{-1})=q(g)\cdot f(n)$ for $g\in L$. The vectors
		$1\otimes(n\Phi(N))$ span $X$ over $k$, so their projections span
		$X_1$. Since $F\ne0$, it follows that $f\ne0$. If $f$ vanished on
		$P_0=[N,T_0]$, then $f(ana^{-1})=f(n)$ for all $a\in T_0$ and $n\in N$.
		Equivariance and $q(T_0)=T$ would give $f(N)\subseteq(\mathcal{E}_1)^T=0$.
		Thus $f(n)\ne0$ for some $n\in P_0$.
		
		We construct a map $f^{L}:L\to\mathcal{E}_1$ satisfying $f^{L}|_N=f$ and
		$f^{L}(gh)=f^{L}(g)+q(g)\cdot f^{L}(h)$. Every element of $q^{-1}(T)$ has a unique
		expression $na$ with $n\in N$ and $a\in T_0$. Define
		$\widetilde f(na)=f(n)$. Then $\widetilde f(nz)=f(n)+\widetilde f(z)$
		for $n\in N$ and $z\in q^{-1}(T)$. Since $q^{-1}(T)$ is normal in $L$,
		we may set
		\[
		f^{L}(g)=\frac{1}{|T_0|}\sum_{a\in T_0}\widetilde f(gag^{-1}).
		\]
		To check the product rule for $f^{L}$, fix $g,h\in L$ and write
		$hah^{-1}=n_a b_a$ with $n_a\in N$ and $b_a\in T_0$. As $a$ ranges
		over $T_0$, so does $b_a$, because $q(b_a)=q(h)q(a)q(h)^{-1}$ and
		$q:T_0\to T$ is an isomorphism. Now
		\[
		\widetilde f(ghah^{-1}g^{-1})
		=f(gn_ag^{-1})+\widetilde f(gb_ag^{-1})
		=q(g)\cdot f(n_a)+\widetilde f(gb_ag^{-1}).
		\]
		Averaging over $a$ and using $f(n_a)=\widetilde f(hah^{-1})$ gives
		$f^{L}(gh)=q(g)\cdot f^{L}(h)+f^{L}(g)$. For $n\in N$ and $a\in T_0$, the
		factorization $nan^{-1}=(nan^{-1}a^{-1})a$ gives
		$\widetilde f(nan^{-1})=f(n)-q(a)\cdot f(n)$. Hence
		$f^{L}(n)=f(n)-e_Tf(n)=f(n)$, since $f(n)\in\mathcal{E}_1$.
		
		Finally, set $\delta(u_g)=f^{L}(g)\rho(u_g)$ and extend $\delta$ by
		$k$-linearity to $A$. The formula for $f^{L}(gh)$, together with the
		conjugation action on $\mathcal{E}$, gives
		\[
		\begin{aligned}
			\delta(u_gu_h)
			&=f^{L}(gh)\rho(u_g)\rho(u_h)\\
			&=f^{L}(g)\rho(u_g)\rho(u_h)+\rho(u_g)f^{L}(h)\rho(u_h)\\
			&=\delta(u_g)\rho(u_h)+\rho(u_g)\delta(u_h).
		\end{aligned}
		\]
		The first equality uses $u_gu_h=\alpha(g,h)u_{gh}$. By bilinearity
		the same product rule holds for all elements of $A$. For $n\in N$,
		$\rho(u_n)=\id_V$, so $\delta(u_n)=f^{L}(n)=f(n)$. The two assertions
		about $\delta(u_n)$ now follow from the corresponding properties of $f$.
	\end{proof}
	
	\subsection{Group transfer}
	
	Let $G_1$ be a subgroup of a finite group $G$, and put $G_1'=[G_1,G_1]$.
	The \emph{group transfer} is a homomorphism
	$\operatorname{Ver}_{G,G_1}:G\to G_1/G_1'$.
	Choose a set $\mathcal T$ of representatives for the left cosets of $G_1$ in $G$.
	For $s\in G$ and $t\in\mathcal T$, write $st=t_s b_{s,t}$ with
	$t_s\in\mathcal T$ and $b_{s,t}\in G_1$. Then
	\[
	\operatorname{Ver}_{G,G_1}(s)=\prod_{t\in\mathcal T}(b_{s,t}G_1').
	\]
	The product is taken in the abelian group $G_1/G_1'$ and is independent
	of the choice of $\mathcal T$ \cite[Theorems~5.1 and~5.2]{Isaacs2008}.
	
	For a homomorphism $\lambda:G_1\to(k,+)$, let $\operatorname{Tr}_{G_1}^G\lambda$
	be the composite of transfer with the map $G_1/G_1'\to(k,+)$ induced by
	$\lambda$. Thus $(\operatorname{Tr}_{G_1}^G\lambda)(s)
	=\sum_{t\in\mathcal T}\lambda(b_{s,t})$. We use transitivity of transfer
	\cite[Theorem~10.8]{Isaacs2008}. For a finite $p$-group $Q$, we have
	$\Phi(Q)=[Q,Q]Q^p$. Every homomorphism $\lambda:Q\to(k,+)$ kills
	commutators and satisfies $\lambda(x^p)=p\lambda(x)=0$ for $x\in Q$.
	Thus $\lambda$ factors through $Q/\Phi(Q)$.
	For an abelian $p$-group $R$, write $\Omega_1(R)=\{r\in R:r^p=1\}$.
	
	\begin{lemma}\label{lem:orbit-sum}
		Let $\partial:A\to A$ be a derivation and let $P\le N$. For $s\in P$, write
		\[
		w(s)=\partial(u_s)u_s^{-1}=\sum_{x\in L}a_x(s)u_x,\qquad a_x(s)\in k.
		\]
		For $g\in L$, put $Q=C_P(g)$ and $g^P=\{tgt^{-1}:t\in P\}$.
		Then $\lambda_g:Q\to(k,+)$, defined by $\lambda_g(b)=a_g(b)$,
		is a group homomorphism, and
		\[
		\sum_{x\in g^P} a_x(s)\rho(u_x)
		=
		(\operatorname{Tr}_{Q}^{P}\lambda_g)(s)\,\rho(u_g)
		\qquad(s\in P).
		\]
	\end{lemma}
	
	\begin{proof}
		Since $P\leq N$ and $\alpha$ is inflated from the normalized cocycle
		$\beta$ on $L/N$, we have $\alpha(s,x)=\alpha(x,s)=1$ for $s\in P$
		and $x\in L$. In particular, $u_{st}=u_su_t$ for $s,t\in P$. Hence
		\[
		\begin{aligned}
			w(st)
			&=\partial(u_su_t)\,u_t^{-1}u_s^{-1}\\
			&=\bigl(\partial(u_s)u_t+u_s\partial(u_t)\bigr)u_t^{-1}u_s^{-1}\\
			&=\partial(u_s)u_s^{-1}+u_s\bigl(\partial(u_t)u_t^{-1}\bigr)u_s^{-1}\\
			&=w(s)+u_s w(t)u_s^{-1}.
		\end{aligned}
		\]
		Also $u_su_xu_s^{-1}=u_{sxs^{-1}}$, so the coefficient of $u_y$ in
		$u_sw(t)u_s^{-1}$ is $a_{s^{-1}ys}(t)$. Comparing coefficients of $u_y$ gives
		\[
		a_y(st)=a_y(s)+a_{s^{-1}ys}(t)\qquad(s,t\in P,\ y\in L).
		\]
		For $b,c\in Q=C_P(g)$ this gives $a_g(bc)=a_g(b)+a_g(c)$, since
		$b^{-1}gb=g$. Thus $\lambda_g$ is a homomorphism.
		
		Choose representatives $\mathcal T$ for the left cosets of $Q$ in $P$.
		As $t$ runs through $\mathcal T$, the elements $g_t=tgt^{-1}$ run through
		$g^P$ exactly once. Fix $s\in P$ and write $st=t_s b_{s,t}$ as above.
		Since $b_{s,t}$ centralizes $g$, we have
		$s g_t s^{-1}=(st)g(st)^{-1}=t_sgt_s^{-1}=g_{t_s}$.
		Apply the coefficient formula with $y=g_{t_s}$ to each of the two
		expressions for $st$. This gives
		\[
		\begin{aligned}
			a_{g_{t_s}}(st)&=a_{g_{t_s}}(s)+a_{g_t}(t),\\
			a_{g_{t_s}}(t_s b_{s,t})&=a_{g_{t_s}}(t_s)+a_g(b_{s,t}).
		\end{aligned}
		\]
		The left sides are equal. Sum the right sides over $t\in\mathcal T$.
		The terms $a_{g_t}(t)$ and $a_{g_{t_s}}(t_s)$ have the same sum, since
		$t\mapsto t_s$ permutes $\mathcal T$, and $a_g(b_{s,t})=\lambda_g(b_{s,t})$.
		Thus $\sum_{x\in g^P}a_x(s)=\sum_{t\in\mathcal T}\lambda_g(b_{s,t})
		=(\operatorname{Tr}_{Q}^P\lambda_g)(s)$.
		Finally, $\rho(u_{g_t})=\rho(u_tu_gu_t^{-1})=\rho(u_g)$ since
		$\rho(u_t)=1_V$. Multiplying the scalar identity by $\rho(u_g)$ gives the formula.
	\end{proof}
	
	Choose one representative from each $P$-conjugacy class in $L$ and denote
	the set of representatives by $\mathcal R_P$.
	For $s\in P$, we have $\rho(\partial(u_s))=\rho(w(s))$ since $\rho(u_s)=1_V$.
	Grouping the terms of $\rho(w(s))$ by these classes and applying
	Lemma~\ref{lem:orbit-sum} to each class gives
	\[
	\rho(\partial(u_s))
	=\sum_{g\in\mathcal R_P}\ \sum_{x\in g^P} a_x(s)\rho(u_x)
	=\sum_{g\in\mathcal R_P}
	(\operatorname{Tr}_{C_P(g)}^P\lambda_g)(s)\,\rho(u_g).
	\]
	
	The next lemma gives a sufficient condition for
	$\operatorname{Tr}_{Q}^P\lambda_g=0$ when $Q=C_P(g)$.
	
	\begin{lemma}\label{lem:elementary-transfer}
		Let $P$ be a finite $p$-group, $Q\le P$, and $R\unlhd P$ elementary abelian.
		If $|R:R\cap Q|>2$, then the image of $\operatorname{Ver}_{P,Q}$
		lies in $\Phi(Q)/[Q,Q]$. In particular, $\operatorname{Tr}_{Q}^P\lambda=0$
		for every homomorphism $\lambda:Q\to(k,+)$.
	\end{lemma}
	
	\begin{proof}
		Put $m=|R:R\cap Q|$. This index is a power of $p$. Thus for odd $p$
		the hypothesis $m>2$ means $R\nleq Q$, and for $p=2$ it means $m\ge4$.
		Since $R$ is normal in $P$, we have $Q\le RQ\le P$. By transitivity
		of transfer \cite[Theorem~10.8]{Isaacs2008} applied to this chain,
		the image of $\operatorname{Ver}_{P,Q}$ is contained in the image of
		$\operatorname{Ver}_{RQ,Q}$. Define $\tau$ to be the composite
		\[
		RQ\xrightarrow{\operatorname{Ver}_{RQ,Q}} Q/Q'
		\xrightarrow{\pi} Q/\Phi(Q),
		\qquad \pi(vQ')=v\Phi(Q).
		\]
		The quotient map $\pi$ is well defined since $Q'=[Q,Q]\le\Phi(Q)$.
		Its kernel is $\Phi(Q)/Q'$, so it suffices to show that
		$\tau(x)=\Phi(Q)$ for every $x\in RQ$.
		
		View $R$ as an $\Fp$-vector space. Choose a complement $R_0$ to $R\cap Q$
		in $R$, so $R=(R\cap Q)\times R_0$ and $|R_0|=m$.
		Then $RQ=R_0Q$ and $R_0\cap Q=1$, so each left coset of $Q$ in
		$RQ$ contains exactly one $t\in R_0$. Fix $y\in R_0$. For each $t\in R_0$,
		write $yt=t_y b_{y,t}$ with $t_y\in R_0$ and $b_{y,t}\in Q$, as in
		the transfer formula. Since $R_0$ is a subgroup, $yt\in R_0$ is itself
		the chosen representative of $ytQ$. Thus $t_y=yt$ and
		$b_{y,t}=t_y^{-1}yt=1$. Hence
		\[
		\tau(y)=\left(\prod_{t\in R_0}b_{y,t}\right)\Phi(Q)=\Phi(Q),
		\]
		the identity of $Q/\Phi(Q)$.
		
		Now fix $u\in Q$. For each $t\in R_0$, write $utu^{-1}=t'f_u(t)$ with
		$t'\in R_0$ and $f_u(t)\in R\cap Q$. The direct product decomposition
		of $R$ gives a homomorphism $R\to R\cap Q$ sending $t'w$ to $w$
		for $t'\in R_0$ and $w\in R\cap Q$. Composing conjugation by $u$ on
		$R_0$ with this map gives $f_u$, so $f_u$ is a homomorphism.
		Multiplying $utu^{-1}=t'f_u(t)$ on the right by $u$ gives $ut=t'(f_u(t)u)$.
		Since $f_u(t)u\in Q$, $t'$ is the chosen representative of $utQ$.
		Thus the element $b_{u,t}$ in the transfer formula is
		$(t')^{-1}ut=f_u(t)u$. As $Q/\Phi(Q)$ is abelian, we obtain
		\[
		\begin{aligned}
			\tau(u)&=\prod_{t\in R_0}\bigl(b_{u,t}\Phi(Q)\bigr)
			=\left(\prod_{t\in R_0} f_u(t)\right)u^m\Phi(Q)\\
			&=f_u\!\left(\prod_{t\in R_0}t\right)u^m\Phi(Q).
		\end{aligned}
		\]
		Here $u^m\in\Phi(Q)$ because $p$ divides $m$ and $u^p\in\Phi(Q)$.
		Since $R_0$ is elementary abelian of order $m>2$, the product of its elements is $1$.
		Thus $\tau(u)=1$. As $R_0$ and $Q$ generate $RQ$, $\tau$ is trivial.
		The assertion about $\lambda$ follows from $\Phi(Q)\leq\ker\lambda$.
	\end{proof}
	
	\begin{lemma}\label{lem:coprime-transfer}
		Let $P$ be a finite $p$-group, let $a\ne1$ be an automorphism of $P$
		of order prime to $p$, and put $Q=C_P(a)$. Then the composite
		\[
		P\xrightarrow{\operatorname{Ver}_{P,Q}}Q/Q'
		\longrightarrow Q/\Phi(Q)
		\]
		is trivial. Thus every $x\in P$ maps to the identity coset $\Phi(Q)$.
	\end{lemma}
	
	\begin{proof}
		Put $R=N_P(Q)$.
		The normalizer condition gives $Q<R$
		\cite[Theorem~1.22]{Isaacs2008}, and coprime action gives
		$R=Q[R,a]$ \cite[Lemma~4.28]{Isaacs2008}.
		Since $a$ fixes $Q$ pointwise and $R$ normalizes $Q$, one has
		$[R,a]\le C_R(Q)$.
		
		Let $\tau:R\to Q/\Phi(Q)$ be the composite of transfer from $R$ to $Q$
		and the quotient map $Q/Q'\to Q/\Phi(Q)$. Choose a set $\mathcal T$
		of representatives for the left cosets of $Q$ in $R$. Since $a$ stabilizes
		both $Q$ and $R=N_P(Q)$, the set $a(\mathcal T)$ also represents these cosets.
		Fix $r\in R$. For each $t\in\mathcal T$, write $rt=t_r b_{r,t}$ with
		$t_r\in\mathcal T$ and $b_{r,t}\in Q$. Applying $a$ gives
		$a(r)a(t)=a(t_r)b_{r,t}$, since $a$ fixes $b_{r,t}$.
		Transfer is independent of the choice of representatives
		\cite[Theorem~5.1]{Isaacs2008}. Computing $\tau(r)$ with $\mathcal T$
		and $\tau(a(r))$ with $a(\mathcal T)$ therefore gives
		\[
		\tau(r)=\prod_{t\in\mathcal T}\bigl(b_{r,t}\Phi(Q)\bigr)=\tau(a(r)).
		\]
		Hence $r^{-1}a(r)\in\ker\tau$ for every $r\in R$, so $[R,a]\le\ker\tau$.
		Every element of $R=Q[R,a]$ induces an inner automorphism on $Q$.
		The transfer formula for the normal subgroup $Q\unlhd R$ therefore
		gives
		\[
		\tau(u)=u^{[R:Q]}\Phi(Q)=\Phi(Q)\qquad(u\in Q).
		\]
		Since $R=Q[R,a]$, the homomorphism $\tau$ is trivial. By transitivity of
		transfer for $Q\le R\le P$ \cite[Theorem~10.8]{Isaacs2008}, we have
		\[
		\operatorname{im}(\operatorname{Ver}_{P,Q})
		\le \operatorname{im}(\operatorname{Ver}_{R,Q})
		\le \Phi(Q)/Q'.
		\]
		Thus the composite from $P$ to $Q/\Phi(Q)$ is trivial.
	\end{proof}
	
	\begin{lemma}[Transfer for groups of class at most two]
		\label{lem:class-two-transfer}
		Let $P$ be a finite $p$-group of nilpotency class at most two,
		let $Q\le P$, and let $a\in\Aut(P)$ fix $Q$ pointwise.
		If
		\[
		P\xrightarrow{\operatorname{Ver}_{P,Q}}Q/Q'\longrightarrow Q/\Phi(Q)
		\]
		is nontrivial, then $a\in O_p(\Aut(P))$.
	\end{lemma}
	
	\begin{proof}
		Choose a homomorphism $\varphi:Q\to(\Fp,+)$ with
		$\operatorname{Tr}_{Q}^P\varphi\ne0$. We use induction on $|P|$.
		The result is immediate if $Q=P$, since then $a=1$.
		
		Suppose that $R_1\le Z(P)$ has order $p$ and $R_1\nleq Q$.
		Then $Q\cap R_1=1$, so $QR_1=Q\times R_1$. Use $R_1$ as the set of
		representatives for the left cosets of $Q$ in $QR_1$. For $q\in Q$ and
		$z,t\in R_1$, we have $(qz)t=(zt)q$. Thus each of the $p$ terms in the
		transfer product is $qQ'$, and
		\[
		\operatorname{Ver}_{QR_1,Q}(qz)=q^pQ'\in\Phi(Q)/Q'.
		\]
		By transitivity of transfer for $Q\le QR_1\le P$
		\cite[Theorem~10.8]{Isaacs2008}, the image of
		$\operatorname{Ver}_{P,Q}$ is contained in the image of $\operatorname{Ver}_{QR_1,Q}$.
		Hence the composite from $P$ to $Q/\Phi(Q)$ is trivial, contrary to the
		hypothesis. Thus every central subgroup of order $p$ lies in $Q$, and
		$R_0=\Omega_1(Z(P))\le Q$. Set $G=C_{\Aut(P)}(R_0)$. Since $R_0$ is
		characteristic in $P$, the subgroup $G$ is normal in $\Aut(P)$, and $a\in G$.
		
		Now suppose that $R_0\cap\ker\varphi$ contains a subgroup $R_1$ of order $p$.
		Since $R_1\le Z(P)$, put $\bar P=P/R_1$, $\bar Q=Q/R_1$, and $\bar x=xR_1$
		for $x\in P$. Since $R_1\le\ker\varphi$, the rule
		$\bar\varphi(\bar q)=\varphi(q)$ defines a homomorphism $\bar Q\to(\Fp,+)$.
		Choose representatives $\mathcal T$ for the left cosets of $Q$ in $P$.
		As $R_1\le Q$, their images represent the left cosets of $\bar Q$ in $\bar P$,
		each exactly once. For $x\in P$ and $t\in\mathcal T$, write $xt=t_x b_{x,t}$
		with $t_x\in\mathcal T$ and $b_{x,t}\in Q$. In $\bar P$ this gives
		$\bar x\bar t=\overline{t_x}\,\overline{b_{x,t}}$. The coset formula gives
		\[
		(\operatorname{Tr}_{\bar Q}^{\bar P}\bar\varphi)(\bar x)
		=\sum_{t\in\mathcal T}\bar\varphi(\overline{b_{x,t}})
		=\sum_{t\in\mathcal T}\varphi(b_{x,t})
		=(\operatorname{Tr}_Q^P\varphi)(x).
		\]
		The right side is nonzero for some $x\in P$. Since $\bar\varphi$ factors
		through $\bar Q/\Phi(\bar Q)$, the transfer from $\bar P$ to $\bar Q$
		followed by projection to $\bar Q/\Phi(\bar Q)$ is nontrivial.
		The automorphism $a$ fixes $R_1\le Q$ pointwise, so it induces an automorphism
		$\bar a$ of $\bar P$ fixing $\bar Q$ pointwise. As $|\bar P|<|P|$ and $\bar P$
		has class at most two, induction gives $\bar a\in O_p(\Aut(\bar P))$.
		
		The action $\theta:G\to\Aut(\bar P)$ has a $p$-group kernel. Indeed, each $\sigma\in\ker\theta$ satisfies $\sigma^p=1$, since
		$\sigma(x)=xz$ with $z\in R_1$ fixed by $\sigma$ for each $x\in P$.
		Thus $\theta^{-1}(O_p(\Aut(\bar P)))$ is a normal $p$-subgroup of $G$
		containing $a$, so $a\in O_p(G)\le O_p(\Aut(P))$.
		
		We may therefore assume that $\varphi|_{R_0}$ is injective.
		Then $Z=Z(P)$ is cyclic and $|R_0|=p$. The restriction of $\varphi$ to
		$R_0$ is therefore an isomorphism onto $(\Fp,+)$. Let $z$ be the unique
		element of $R_0$ with $\varphi(z)=1$. Then $R_0=\langle z\rangle$.
		If $P$ is abelian, then $P=Z(P)$ is cyclic. Write $|a|=p^e m$ with
		$p\nmid m$. The automorphism $a^{p^e}$ has order $m$ and fixes
		$\Omega_1(P)=R_0\le Q$ pointwise, so $a^{p^e}=1$ by
		\cite[Corollary~4.35]{Isaacs2008}. Thus $a$ is a $p$-element of the abelian
		group $\Aut(P)$ and hence belongs to $O_p(\Aut(P))$.
		We may assume that $P$ is nonabelian and $Q<P$.
		Every nontrivial subgroup of the cyclic $p$-group $Z$ contains
		$R_0=\langle z\rangle$. Since $\varphi(z)=1$, we have $Z\cap\ker\varphi=1$.
		Thus $\varphi$ embeds $Q\cap Z$ in $(\Fp,+)$, so $|Q\cap Z|\le p$.
		As $R_0\le Q\cap Z$ has order $p$, we obtain $Q\cap Z=R_0$.
		Also $Q'\le P'\le Z$ because $P$ has class at most two, and
		$Q'\le\ker\varphi$ because $(\Fp,+)$ is abelian. Hence $Q'=1$.
		
		Put $R=QZ$, $s=[P:R]$, and $m=\exp(P/Z)$.
		The group $R$ is abelian and normal.
		Write $[g,h]=g^{-1}h^{-1}gh$ and $a^b=b^{-1}ab$. Since $P'\le Z=Z(P)$,
		all commutators are central, so the commutator identities
		\cite[Lemma~2.11(i)--(ii)]{Craven2008} give
		\[
		\begin{aligned}
			[g,hk]&=[g,k][g,h]^k=[g,h][g,k],\\
			[gu,h]&=[g,h]^u[u,h]=[g,h][u,h].
		\end{aligned}
		\]
		For $z,w\in Z$ we have $[gz,hw]=[g,h]$, so $\ell_g(hZ)=[g,h]$
		depends only on $gZ$ and $hZ$. The first identity shows that
		$\ell_g:P/Z\to Z$ is a homomorphism. The second gives
		$\ell_{gu}(hZ)=\ell_g(hZ)\ell_u(hZ)$ for every $hZ\in P/Z$.
		Thus $gZ\mapsto\ell_g$ is a homomorphism from $P/Z$ to $\Hom(P/Z,Z)$.
		This map is injective, since $\ell_g=1$ if and only if $g$
		commutes with every element of $P$, that is, if and only if $g\in Z$.
		Hence $\ell_g$ has the same order as $gZ$.
		
		For $r\in R$, the homomorphism $\ell_r$ is trivial on $R/Z$
		because $R$ is abelian. Hence the restriction of
		$P/Z\hookrightarrow\Hom(P/Z,Z)$ to $R/Z$ induces an injection
		$R/Z\hookrightarrow\Hom(P/R,Z)$.
		As $P/R$ is abelian and $Z$ is cyclic, each cyclic factor of $P/R$
		has at most its order many homomorphisms to $Z$. Thus $|R:Z|\le [P:R]=s$.
		
		Choose $gZ$ of order $m$ in the finite abelian group $P/Z$.
		The image $\ell_g(P/Z)$ is a subgroup of the cyclic group $Z$.
		Choose $h\in P$ such that $[g,h]=\ell_g(hZ)$ generates this image.
		Every value of $\ell_g$ is then a power of $[g,h]$, so for each
		$n\ge1$ we have $\ell_g^n=1$ if and only if $[g,h]^n=1$.
		Hence $[g,h]$ has the same order as $\ell_g$, namely $m$.
		If $(gZ)^i(hZ)^j=Z$, then $g^ih^j\in Z$. Since commutators in $P$ are central,
		\[
		1=[g^ih^j,h]=[g,h]^i,\qquad
		1=[g^ih^j,g]=[h,g]^j=[g,h]^{-j}.
		\]
		Thus $m$ divides both $i$ and $j$. As $m=\exp(P/Z)$, it follows that
		$gZ$ and $hZ$ generate a subgroup isomorphic to $C_m\times C_m$.
		Hence
		\begin{equation}\label{eq:transfer-index}
			s\ge\sqrt{|P:Z|}\ge m.
		\end{equation}
		If $s=m$, then \eqref{eq:transfer-index} gives $|P:Z|=m^2$, so
		$P/Z=\langle gZ,hZ\rangle\cong C_m\times C_m$.
		
		Set $\psi=\operatorname{Tr}_{Q}^R\varphi$.
		Since $R$ is abelian and $[R:Q]=|Z:R_0|$, transfer sends a generator
		$r$ of $Z$ to $r^{[R:Q]}$, a generator of $R_0$.
		Consequently
		\begin{equation}\label{eq:transfer-center-kernel}
			\ker(\psi|_Z)=Z^p.
		\end{equation}
		We compute $(\operatorname{Tr}_R^P\psi)(x)$ by the coset formula.
		If $x\notin R$, every orbit of $\langle x\rangle$ on $P/R$ has
		length $d=|xR|$, with $p\mid d$. By the transfer-evaluation lemma
		\cite[Lemma~5.5(b)]{Isaacs2008},
		the orbit containing $tR$ contributes $t^{-1}x^dt$ to the transfer product.
		As $[x,t]\in Z$, we have $t^{-1}x^dt=(x[x,t])^d=x^d[x,t]^d$.
		Since $[x,t]^d\in Z^p$,
		we have $\psi(t^{-1}x^dt)=\psi(x^d)$ by \eqref{eq:transfer-center-kernel}.
		There are $s/d$ orbits, so $(\operatorname{Tr}_R^P\psi)(x)=\psi(x^s)$.
		
		If $x\in R$, every coset is fixed, and $t^{-1}xt=x[x,t]$. Hence
		\[
		(\operatorname{Tr}_R^P\psi)(x)-\psi(x^s)
		=\sum_{tR\in P/R}\psi([x,t]).
		\]
		The image of $xZ$ under $R/Z\hookrightarrow\Hom(P/R,Z)$ is the homomorphism
		$tR\mapsto[x,t]$. Composing it with $\psi|_Z$ gives the homomorphism
		$tR\mapsto\psi([x,t])$ from $P/R$ to $(\Fp,+)$.
		If it is nonzero, each element of $\Fp$ occurs $s/p$ times in the sum.
		For odd $p$ the elements of $\Fp$ sum to zero. For $p=2$ and $s>2$,
		the multiplicity $s/2$ is even. Thus only $(p,s)=(2,2)$ requires a separate check.
		
		Suppose that $(p,s)=(2,2)$.
		By \eqref{eq:transfer-index}, $P/Z\cong C_2^2$.
		Thus $u^2\in Z$ for every $u\in P$. Since commutators are central,
		$[u,v]^2=[u^2,v]=1$ for all $u,v\in P$. Hence every commutator lies in
		$\langle z\rangle$, the unique subgroup of order two in the cyclic group $Z$.
		If $|Z|\ge4$, then $Z^2$ has even order and contains $\langle z\rangle$.
		Thus $P'\le Z^2=\ker(\psi|_Z)$ by \eqref{eq:transfer-center-kernel},
		so every term $\psi([x,t])$ in the sum is zero.
		If $|Z|=2$, then $Z=\langle z\rangle\le Q$, so $R=QZ=Q$.
		Since $|P:Z|=4$ and $[P:R]=2$, we have $|P|=8$ and $|Q|=4$.
		If $Q=\langle u\rangle$ were cyclic, its unique element of order two
		would be $u^2=z$. This would give $\varphi(z)=2\varphi(u)=0$ in $\mathbb F_2$,
		contrary to $\varphi(z)=1$. Thus $Q\cong C_2\times C_2$.
		It follows that $P\cong D_8$, and $\Aut(P)$ is a
		$2$-group.  The conclusion already holds in this exceptional case.
		
		In all remaining cases the calculation gives
		\begin{equation}\label{eq:transfer-power}
			(\operatorname{Tr}_R^P\psi)(x)=\psi(x^s)\qquad(x\in P).
		\end{equation}
		By transitivity,
		$(\operatorname{Tr}_{Q}^P\varphi)(x)=\psi(x^s)$.
		If $s>m$, then $p\mid s/m$ and $x^s=(x^m)^{s/m}\in Z^p$ for every $x\in P$.
		Since $Z^p=\ker(\psi|_Z)$ by \eqref{eq:transfer-center-kernel}, this gives
		\[
		(\operatorname{Tr}_{Q}^P\varphi)(x)=\psi(x^s)=0\qquad(x\in P),
		\]
		contrary to the choice of $\varphi$ with $\operatorname{Tr}_{Q}^P\varphi\ne0$.
		Thus $s=m$, and some $x$ satisfies $x^m\notin Z^p$.
		Since $Z$ is cyclic, this means $Z=\langle x^m\rangle$.
		The element $xZ$ has order $m$: a smaller order would put $x^m$
		in $Z^p$.
		The injection $P/Z\hookrightarrow\Hom(P/Z,Z)$ sends $xZ$ to $\ell_x$,
		so $\ell_x$ has order $m$. Choose $y\in P$ so that
		$[x,y]=\ell_x(yZ)$ generates the cyclic subgroup $\ell_x(P/Z)$ of $Z$.
		Since $\ell_x^n=1$ exactly when $[x,y]^n=1$, we obtain
		\begin{equation}\label{eq:transfer-commutator-order}
			|[x,y]|=|\ell_x|=|xZ|=m.
		\end{equation}
		The equality case of \eqref{eq:transfer-index} now gives
		\[
		P=\langle x,y\rangle,\qquad P/Z\cong C_m^2,\qquad
		Z=\langle x^m\rangle\le\Phi(P),\qquad
		P/\Phi(P)\cong\Fp^2.
		\]
		
		Since $m=\exp(P/Z)$, we have $g^m\in Z$ for every $g\in P$.
		As $[h,g]\in Z$, the power formula
		\cite[Lemma~3.3]{BlackburnNeumannVenkataraman2007} gives
		\[
		(gh)^m=g^mh^m[h,g]^{\binom{m}{2}}.
		\]
		Now $[h,g]^m=[h,g^m]=1$. For odd $p$, the divisibility
		$m\mid\binom{m}{2}$ makes the last factor trivial. For $p=2$, its square
		is trivial, and $m\mid |Z|$ by \eqref{eq:transfer-commutator-order}.
		Since the case $m=|Z|=2$ was settled above, the cyclic $2$-group $Z$
		has order at least four, and every element of order at most two lies in $Z^2$.
		Thus in both cases $(gh)^mZ^p=(g^mZ^p)(h^mZ^p)$, so the map
		\[
		\eta:P\longrightarrow Z/Z^p,\qquad g\longmapsto g^mZ^p
		\]
		is a homomorphism. It is nonzero since $x^m\notin Z^p$.
		Each $\alpha\in G=C_{\Aut(P)}(R_0)$ acts on the cyclic group $Z$ by
		$u\mapsto u^k$ for some integer $k$ \cite[2.2.2]{KurzweilStellmacher2004}.
		Since $\alpha$ fixes the generator $z$ of $R_0$ and $|z|=p$,
		we have $z^k=z$ and hence $k\equiv1\pmod p$. Thus
		$\alpha(u)Z^p=u^kZ^p=uZ^p$ for every $u\in Z$.
		For $\alpha\in G$ and $g\in P$ we therefore have
		$\eta(\alpha(g))=\alpha(\eta(g))=\eta(g)$.
		
		Let $A\le G$ have order coprime to $p$. By the coprime action theorem
		\cite[Theorem~6.2]{Craven2008}, every $A$-fixed element in the image of an
		$A$-equivariant homomorphism from $P$ has an $A$-fixed preimage.
		Applying this to $\eta$, we may choose $x$ with $x^m\notin Z^p$ fixed by $A$.
		Then $A$ fixes $Z=\langle x^m\rangle$ pointwise. The homomorphism
		$P\to Z$, $t\mapsto[x,t]$, is now $A$-equivariant, since
		$[x,\alpha(t)]=\alpha([x,t])$ for $\alpha\in A$. Applying the same theorem
		to a generator of its cyclic image, we may choose $y$ in
		\eqref{eq:transfer-commutator-order} fixed by $A$ as well.
		Since $P=\langle x,y\rangle$, we have $A=1$.
		
		By Cauchy's theorem, $G$ is a $p$-group. Normality in $\Aut(P)$ now gives
		$a\in G\le O_p(\Aut(P))$.
	\end{proof}
	
	\subsection{The critical subgroup and the proof of the criterion}
	
	For a nontrivial action of $g$ on an elementary abelian $p$-group $Y$,
	the bound $|Y:C_Y(g)|>2$ can fail only at $p=2$ and index two.
	The following lemma excludes this case for the faithful action used below.
	
	\begin{lemma}\label{lem:no-transvection}
		Suppose $p=2$, $H$ is solvable, $O_2(H)=1$, and $k^\beta H$ has
		just one simple module. If $M$ is a faithful finite-dimensional
		$\mathbb F_2H$-module, then $\rank_M(s-1)\ne1$ for every $s\in H$.
	\end{lemma}
	
	\begin{proof}
		Write $(w_g)_{g\in H}$ for the standard $k$-basis of the twisted group
		algebra $k^\beta H$, with multiplication $w_a w_b=\beta(a,b)w_{ab}$ for
		$a,b\in H$. An element $g\in H$ is called $\beta$-regular if
		$\beta(g,h)=\beta(h,g)$ for every $h\in C_H(g)$. For such $h$ we have
		$gh=hg$, so $w_g w_h=\beta(g,h)w_{gh}$ and $w_h w_g=\beta(h,g)w_{hg}$
		are multiples of the same basis element. Thus $g$ is $\beta$-regular
		exactly when $w_g$ commutes with $w_h$ for every $h\in C_H(g)$.
		
		By Conlon's counting formula \cite[p.~156]{Conlon1964}, the number
		of isomorphism classes of simple $k^\beta H$-modules equals the number
		of $\beta$-regular $H$-conjugacy classes of elements of odd order.
		The identity already gives one such class. Since $k^\beta H$ has only
		one simple module up to isomorphism, the identity class is the only such class.
		
		Suppose, for a contradiction, that $\rank_M(s-1)=1$ for some $s\in H$.
		We will construct a nonidentity $\beta$-regular element $z$ of odd order
		in the Fitting subgroup $F=F(H)$, giving a second such class.
		Put $C=C_M(s)=\ker(s-1)$. By the rank assumption, $M/C$ has dimension
		one over $\mathbb F_2$. Its only nonzero coset must be fixed by $s$.
		Thus $sv+C=v+C$ for every $v\in M$, so $(s-1)M\subseteq C=\ker(s-1)$.
		It follows that $s^2-1=(s-1)^2=0$ on $M$. Since $s$ acts
		nontrivially and $M$ is faithful, $s$ has order two.
		
		The Sylow $2$-subgroup of $F$ is characteristic in $F$ and hence
		normal in $H$. Since $O_2(H)=1$, $F$ has odd order. Solvability gives
		$C_H(F)\le F$ \cite[Theorem~9.8]{Isaacs2008}, so the involution $s$ cannot centralize $F$.
		Choose $a\in F$ with $t=a^{-1}sa\ne s$, and put $z=st$.
		Since $F$ is normal, $z=(sa^{-1}s^{-1})a$ lies in $F$. Thus $z$ has
		odd order and is nontrivial because $s$ and $t$ are distinct involutions.
		
		Let $n=\dim M$. Since $t$ is conjugate to $s$, both $s-1$ and $t-1$
		have rank one on $M$. By rank-nullity, $\dim C_M(s)=\dim C_M(t)=n-1$.
		The dimension formula gives
		$\dim(C_M(s)\cap C_M(t))\ge (n-1)+(n-1)-n=n-2$.
		Every vector in this intersection is fixed by $z=st$, so
		$W=C_M(z)$ has dimension at least $n-2$. The odd-order group $\langle z\rangle$
		acts on the elementary abelian $2$-group $(M,+)$, so coprime action
		\cite[Theorem~4.34]{Isaacs2008} gives $M=W\oplus U$, where
		$U=[M,\langle z\rangle]=(z-1)M$. Thus $\dim U=n-\dim W\le 2$. Faithfulness gives
		$U\ne0$. If $\dim U=1$, then $z$ would fix the only nonzero vector
		of $U$ and hence all of $M$. Thus $\dim U=2$. The action of $z$ on the three
		nonzero vectors of $U$ is a nonidentity permutation of odd order,
		so it is a $3$-cycle. Since $z$ fixes $W$ and $M$ is faithful, $|z|=3$.
		
		Now $\dim W=n-2$. The subspace $C_M(s)\cap C_M(t)$ lies in $W$
		and has dimension at least $n-2$, so it equals $W$. Thus $s$ and $t$
		fix $W$ pointwise. Since $s$ and $t$ normalize $\langle z\rangle$,
		they preserve $U=[M,\langle z\rangle]$.
		Take $h\in C_H(z)$. Then $h$ preserves both $U$ and $W$. On the three
		nonzero vectors of $U$, the only permutations commuting with the $3$-cycle
		$z|_U$ are the powers of $z|_U$.
		Hence $h|_U=z^i|_U$ for some $i$. Put $h_0=z^{-i}h$. Then $h_0$
		fixes $U$ pointwise and preserves $W$. For $r\in\{s,t\}$, $u\in U$
		and $w\in W$, we have $h_0r(u+w)=ru+h_0w=rh_0(u+w)$, since $r$
		preserves $U$ and fixes $W$ pointwise. Thus $h_0r$ and $rh_0$ induce
		the same linear map on $M=U\oplus W$. Since $M$ is a faithful $H$-module,
		$h_0r=rh_0$ in $H$.
		
		It remains to show that $w_z$ commutes with $w_h$. For $r=s,t$ we have
		$h_0rh_0^{-1}=r$, so $w_{h_0}w_rw_{h_0}^{-1}=\lambda w_r$ for some
		$\lambda\in k^\times$. Squaring gives $w_r^2=\lambda^2w_r^2$,
		because $r^2=1$ makes $w_r^2$ a nonzero scalar. In characteristic two,
		$\lambda^2=1$ forces $\lambda=1$. Thus $w_{h_0}$ commutes with
		$w_sw_t$ and hence with its scalar multiple $w_z$. Since $h=z^ih_0$,
		$w_h$ is a scalar multiple of $w_z^i w_{h_0}$ and also commutes with
		$w_z$. We have therefore proved, for every $h\in C_H(z)$, that
		\[
		w_z w_h=w_h w_z=\lambda_h w_{zh},\qquad
		\lambda_h=\beta(z,h)=\beta(h,z).
		\]
		Thus $1$ and $z$ represent distinct $\beta$-regular conjugacy
		classes of elements of odd order, contrary to the hypothesis that $k^\beta H$ has
		only one simple module.
	\end{proof}
	
	\begin{proof}[Proof of Proposition~\ref{prop:derivation}]
		\emph{(i) implies (ii).}
		If $L$ is a $p$-group, then $H=1$, so our choice of cocycle gives
		$A=kL$. Its simple module affords the augmentation $\varepsilon:kL\to k$.
		Let $\delta:kL\to k$ be an $\varepsilon$-derivation. Since
		$\varepsilon(g)=1$ for $g\in L$, we have $\delta(gh)=\delta(g)+\delta(h)$.
		Define $\partial:kL\to kL$ on the group basis by $\partial(g)=\delta(g)g$
		and extend linearly \cite[p.~341]{FarkasGreenMarcos2000}. Then
		\[
		\partial(gh)=(\delta(g)+\delta(h))gh
		=\partial(g)h+g\partial(h).
		\]
		Thus $\partial$ is a derivation and $\varepsilon\circ\partial=\delta$.
		Lemma~\ref{lem:lift-derivation} gives~(ii).
		
		\emph{(ii) implies (i).}
		Assume~(ii). Let $V$ be the unique simple $A$-module and let
		$\rho:A\to\mathcal E=\End_k(V)$ be its representation.
		By Lemma~\ref{lem:one-simple-reduction}, $H=L/N$ is solvable,
		$O_p(H)=1$, and $C_L(N)\le N$. If $N=1$, the last inclusion gives
		$L=1$. We may therefore assume $N\ne1$. Suppose, towards a
		contradiction, that $H\ne1$.
		
		Choose a Thompson critical subgroup $Q$ of $N$. Then $Q$ is
		characteristic in $N$, $[Q,N]\le Z(Q)$, and $C_N(Q)=Z(Q)$
		\cite[Theorem~6.5]{Craven2008}. In particular, $Y=\Omega_1(Z(Q))$ is a
		normal elementary abelian subgroup of $L$. Choose a composition series
		$0=Y_0<Y_1<\cdots<Y_r=Y$ of $\Fp L$-modules and put
		$S_i=Y_i/Y_{i-1}$. Since $S_i$ is simple and $N\unlhd L$, the restriction
		$S_i\!\downarrow_N$ is semisimple by Clifford's theorem
		\cite[Theorem~5.3.1]{Webb2016}. As $N$ is a $p$-group, the trivial
		module is the only simple $\Fp N$-module up to isomorphism
		\cite[Proposition~6.2.1]{Webb2016}. Hence $N$ acts trivially on each $S_i$,
		and the action of $L$ on $M=\bigoplus_i S_i$ factors through $H=L/N$.
		For $g\in L$, write $\bar g=gN$. The projection
		$Y_i\to S_i$ sends $C_{Y_i}(g)$ into $C_{S_i}(\bar g)$ with kernel
		$C_{Y_{i-1}}(g)$. Summing dimensions and using $\dim Y=\dim M$ gives
		\begin{equation}\label{eq:rank-composition}
			\rank_Y(g-1)\ge \rank_M(\bar g-1).
		\end{equation}
		
		If the action of $H$ on $M$ has nontrivial kernel, choose a minimal
		normal subgroup $T$ of $H$ contained in this kernel. If the action is
		faithful, take any minimal normal subgroup $T$ of $H$.
		Since $H$ is solvable and $O_p(H)=1$, $T$ is an elementary abelian
		$p'$-group. By Schur--Zassenhaus \cite[Theorem~3.8]{Isaacs2008}, the inverse
		image of $T$ in $L$ has the form $N\rtimes T_0$. Put $P_0=[N,T_0]$
		and let $e_T=|T|^{-1}\sum_{t\in T}t$ act on $\mathcal E$ by averaging
		the conjugation action of $T$. Lemma~\ref{lem:target-derivation}
		supplies a $\rho$-derivation $\delta$ and an element $n_0\in P_0$ with
		$0\ne\delta(u_{n_0})\in(1-e_T)\mathcal E$. By~(ii) and
		Lemma~\ref{lem:lift-derivation}, $\delta=\rho\circ\partial_0$ for some
		derivation $\partial_0$ of $A$. Since $\mathcal E^T=e_T\mathcal E$ and
		$e_T$ is idempotent, it suffices for a contradiction to prove that every
		derivation $\partial$ of $A$ satisfies
		\begin{equation}\label{eq:vanish-on-commutators}
			\rho(\partial(u_n))\in\mathcal E^T \qquad(n\in P_0).
		\end{equation}
		
		Fix a derivation $\partial$ of $A$. For $s\in N$, write
		$\partial(u_s)u_s^{-1}=\sum_{x\in L}a_x(s)u_x$. Since $\rho(u_s)=1_V$,
		applying $\rho$ gives
		\begin{equation}\label{eq:derivation-expansion}
			\rho(\partial(u_s))=\sum_{x\in L}a_x(s)\rho(u_x).
		\end{equation}
		For $g\in L$, the map $\lambda_g:C_N(g)\to(k,+)$ given by
		$\lambda_g(b)=a_g(b)$ is a homomorphism by Lemma~\ref{lem:orbit-sum}.
		For $P\le N$ and $s\in P$, the same lemma gives the contribution
		of the $P$-conjugacy class $g^P$ to \eqref{eq:derivation-expansion} as
		\begin{equation}\label{eq:orbit-contribution}
			\sum_{x\in g^P}a_x(s)\rho(u_x)
			=\bigl(\operatorname{Tr}_{C_P(g)}^P
			(\lambda_g|_{C_P(g)})\bigr)(s)\,\rho(u_g).
		\end{equation}
		
		\smallskip\noindent\textbf{\boldmath Case 1. $H$ acts faithfully on $M$.}
		For $g\in L\setminus N$, faithfulness gives $\rank_M(\bar g-1)\ge1$.
		If $p=2$, Lemma~\ref{lem:no-transvection} excludes rank one, since
		$H$ is solvable, $O_2(H)=1$, and $k^\beta H$ has a unique simple
		module. Thus \eqref{eq:rank-composition} gives
		\[
		|Y:C_Y(g)|=p^{\rank_Y(g-1)}
		\ge \begin{cases}4 & \text{if }p=2,\\ p & \text{if }p\text{ is odd}.\end{cases}
		\]
		In either case $|Y:Y\cap C_N(g)|>2$. Applying
		Lemma~\ref{lem:elementary-transfer} to $Y\unlhd N$ and $C_N(g)\le N$
		gives $\operatorname{Tr}_{C_N(g)}^N\lambda_g=0$. By
		\eqref{eq:orbit-contribution} with $P=N$, every $N$-conjugacy class
		in $L\setminus N$ contributes zero to \eqref{eq:derivation-expansion}.
		The remaining terms have $x\in N$ and hence $\rho(u_x)=1_V$. Thus
		$\rho(\partial(u_n))\in k1_V\subseteq\mathcal E^T$ for every $n\in N$.
		This proves \eqref{eq:vanish-on-commutators} in Case~1.
		
		\smallskip\noindent\textbf{\boldmath Case 2. The action of $H$ on $M$ has nontrivial kernel.}
		By the choice of $T$, the complement $T_0$ acts trivially on every
		$Y_i/Y_{i-1}$. In a basis adapted to this series, its image on $Y$
		is unitriangular and hence a $p$-group. Since $T_0$ is a $p'$-group,
		it centralizes $Y=\Omega_1(Z(Q))$. Coprime action on the abelian
		$p$-group $Z(Q)$ now gives $[T_0,Z(Q)]=1$
		\cite[Corollary~4.35]{Isaacs2008}.
		
		Put $K_0=C_N(Z(Q))$. Since $T_0$ centralizes $Z(Q)$, $n$ and
		$ana^{-1}$ induce the same automorphism of $Z(Q)$ for $n\in N$ and
		$a\in T_0$. Thus $P_0=[N,T_0]\le K_0$. For each $x\in K_0$, let
		$c_x:Q\to Q$ denote conjugation by $x$. Then $c_x$ fixes $Z(Q)$
		pointwise and induces the identity on $Q/Z(Q)$ since $[Q,N]\le Z(Q)$.
		Hence $c_x(q)=xqx^{-1}=q\,d_x(q)$ for $q\in Q$, where
		$d_x:Q\to Z(Q)$ is a homomorphism trivial on $Z(Q)$.
		For $x,y\in K_0$ and $q\in Q$, both $d_x(q)$ and $d_y(q)$ belong
		to $Z(Q)$, which $K_0$ fixes pointwise. Hence
		\[
		\begin{aligned}
			(c_xc_y)(q)&=c_x\bigl(q\,d_y(q)\bigr)=q\,d_x(q)d_y(q),\\
			(c_yc_x)(q)&=c_y\bigl(q\,d_x(q)\bigr)=q\,d_y(q)d_x(q).
		\end{aligned}
		\]
		The two products agree because $d_x(q)$ and $d_y(q)$ commute.
		Thus conjugation by $[x,y]$ induces $[c_x,c_y]=\mathrm{id}_Q$,
		so every commutator in $K_0$ centralizes $Q$. It follows that
		\begin{equation}\label{eq:central-commutators}
			[K_0,K_0]\le C_N(Q)=Z(Q)\le Z(K_0).
		\end{equation}
		In particular, $P_0$ has nilpotency class at most two.
		
		By \cite[Lemma~4.1]{Isaacs2008}, $N$ normalizes $P_0=[N,T_0]$.
		Conjugation by $L$ sends $T_0$ to another complement in the inverse image of $T$,
		and all such complements are $N$-conjugate \cite[Theorem~3.12]{Isaacs2008}.
		Hence $P_0\unlhd L$. By coprime action,
		$N=P_0C_N(T_0)$ \cite[Lemma~4.28]{Isaacs2008}. An element of $T_0$
		centralizing $P_0$ would then centralize $N$, so it lies in
		$C_L(N)\cap T_0=1$. Thus $T_0$ acts faithfully on $P_0$.
		
		Let $G_0$ be the image of $L$ in $\Aut(P_0)$. The image of $N$ is
		a normal $p$-subgroup of $G_0$, so the action induces a homomorphism
		$H\to G_0/O_p(G_0)$. Its restriction to $T$ is injective because
		$T_0$ acts faithfully and has order prime to $p$. Let $g\in L$ and
		write $\bar g=gN$ for its image in $H=L/N$. If $\bar g\notin C_H(T)$,
		then the automorphism $u\mapsto gug^{-1}$ of $P_0$ lies outside
		$O_p(\Aut(P_0))$. Indeed, $G_0\cap O_p(\Aut(P_0))\le O_p(G_0)$.
		If $g$ induced an automorphism in this intersection, then $\bar g$
		would lie in the kernel of $H\to G_0/O_p(G_0)$. For every $t\in T$,
		the commutator $[\bar g,t]$ would then lie both in that kernel and
		in $T$, forcing $[\bar g,t]=1$.
		
		Fix $n\in P_0$ and take $s=n$ in \eqref{eq:derivation-expansion}.
		First collect the terms indexed by elements $g$ with nontrivial $p'$-part $a$ into
		$N$-conjugacy classes. Since $C_L(N)\le N$, $a$ induces a nontrivial
		automorphism of $N$ of order prime to $p$. As $C_N(g)\le C_N(a)$,
		transitivity gives
		\[
		\operatorname{Tr}_{C_N(g)}^N\lambda_g
		=\operatorname{Tr}_{C_N(a)}^N
		\bigl(\operatorname{Tr}_{C_N(g)}^{C_N(a)}\lambda_g\bigr)=0
		\]
		by Lemma~\ref{lem:coprime-transfer}. By \eqref{eq:orbit-contribution}
		with $P=N$, each of these classes contributes zero.
		
		Group the remaining terms, indexed by $p$-elements $g$, into
		$P_0$-conjugacy classes. Suppose $\bar g\in C_H(T)$ and take $t\in T$.
		In $k^\beta H$ we have
		$w_tw_{\bar g}w_t^{-1}=c\,w_{\bar g}$ for some $c\in k^\times$.
		Since $w_{\bar g}^{|\bar g|}$ is a nonzero scalar,
		$c^{|\bar g|}=1$. The order of $\bar g$ is a power of $p$, so $c=1$
		in characteristic $p$. Thus $w_t$ commutes with $w_{\bar g}$, which
		acts on $V$ as $\rho(u_g)$. Hence $t\cdot\rho(u_g)=\rho(u_g)$.
		Taking $P=P_0$ and $s=n$ in \eqref{eq:orbit-contribution} gives
		\[
		\sum_{x\in g^{P_0}}a_x(n)\rho(u_x)
		=\bigl(\operatorname{Tr}_{C_{P_0}(g)}^{P_0}
		(\lambda_g|_{C_{P_0}(g)})\bigr)(n)\,\rho(u_g)
		\in\mathcal E^T.
		\]
		Indeed, the transfer coefficient lies in $k$, so the right-hand side
		is a scalar multiple of the $T$-fixed operator $\rho(u_g)$.
		
		If $\bar g\notin C_H(T)$, conjugation by $g$ fixes $C_{P_0}(g)$
		pointwise and lies outside $O_p(\Aut(P_0))$. Since $P_0\le K_0$
		has class at most two by \eqref{eq:central-commutators},
		Lemma~\ref{lem:class-two-transfer} implies that transfer
		from $P_0$ to $C_{P_0}(g)$ is trivial modulo $\Phi(C_{P_0}(g))$.
		Every homomorphism from $C_{P_0}(g)$ to $(k,+)$ factors through this
		Frattini quotient, so the coefficient of $\rho(u_g)$ in
		\eqref{eq:orbit-contribution} with $P=P_0$ is zero. All contributions
		to \eqref{eq:derivation-expansion} therefore lie in $\mathcal E^T$, proving
		\eqref{eq:vanish-on-commutators} in Case~2.
		
		In either case, applying \eqref{eq:vanish-on-commutators} to
		$\partial_0$ and $n_0$ gives
		$\delta(u_{n_0})\in\mathcal E^T\cap(1-e_T)\mathcal E=0$,
		contrary to the choice of $n_0$. Hence $H=1$, and $L=N$ is a $p$-group.
	\end{proof}
	
	\begin{corollary}\label{cor:twisted}
		Let $O_{p'}(L)=1$. If $k^\alpha L$ is Morita equivalent to $kP$ for
		a finite $p$-group $P$, then $L$ is a $p$-group.
	\end{corollary}
	
	\begin{proof}
		Let $B$ be the basic algebra of $k^\alpha L$. Since $kP$ is local
		and basic, Morita equivalence implies that $k^\alpha L$ has a unique
		simple module up to isomorphism and that $B\cong kP$
		\cite[Corollary~4.9.6]{Linckelmann2018a}.
		By Proposition~\ref{prop:derivation} applied to the $p$-group $P$,
		$kP$ satisfies condition~\textup{(ii)}, and hence so does $B$.
		Since $O_{p'}(L)=1$, the reverse implication applied to $L$ shows
		that $L$ is a $p$-group.
	\end{proof}
	
	\section{Proof of Theorem~\ref*{thm:nilpotency}}\label{sec:integral}
	
	We now combine Theorem~\ref{thm:solvable} with Corollary~\ref{cor:twisted}.
	For a normalized cocycle $\alpha$ with values in $\cO^\times$, the notation $\cO^\alpha L$ means the free $\cO$-module with basis $\{u_g:g\in L\}$ and multiplication $u_gu_h=\alpha(g,h)u_{gh}$.
	
	\begin{lemma}\label{lem:fong}
		Let $A$ be a block algebra of a finite solvable group over $\cO$.
		There are a finite solvable group $L$ with $O_{p'}(L)=1$ and a normalized cocycle
		\[
		\alpha\in Z^2(L,\mu_m(\cO)),\qquad (m,p)=1,
		\]
		such that $A\sim_{\Mor}\cO^\alpha L$.
		Here $\mu_m(\cO)$ denotes the group of $m$th roots of unity in $\cO$.
	\end{lemma}
	
	\begin{proof}
		Choose a solvable group $X$ of least order with a block $c$ such that
		$A\sim_{\Mor}\cO Xc$. Put $N=O_{p'}(X)$ and $L=X/N$, so $O_{p'}(L)=1$.
		Let $f$ be a block of $\cO N$ covered by $c$. By Fong--Reynolds
		reduction \cite[Theorem~6.8.3]{Linckelmann2018b} and the minimality of $|X|$, $f$ is
		$X$-invariant. Fong's theorem \cite[Theorem~10.20]{Navarro1998} makes $c$ the
		unique block of $\cO X$ covering $f$, so $c=f$ in $\cO X$.
		
		Since $N$ is a $p'$-group, $f$ has trivial defect group and is nilpotent.
		The model group in \cite[Theorem~1.1(b)]{CoconetMarcusTodea2020} is therefore $X/N$, and
		the theorem gives an $L$-graded basic Morita equivalence between
		$\cO Xf$ and $\cO^\alpha L$.
		With $m=|L|_{p'}$, the Clifford class in $H^2(L,k^\times)$ has a
		normalized representative with values in $\mu_m(k)$
		\cite[Proposition~1.2.9(ii)]{Linckelmann2018a}. By \cite[p.~6 of the arXiv version]{CoconetMarcusTodea2020}, we may take
		$\alpha$ to be its multiplicative lift to $\mu_m(\cO)$.
	\end{proof}
	
	\begin{proof}[Proof of Theorem~\ref{thm:nilpotency}]
		As in Section~\ref{sec:preliminaries}, $b$ denotes both the given block of $kG$ and its unique lift to $\cO G$.
		Since $kGb\sim_{\Mor}kP$, we have $l(b)=1$.
		By Theorem~\ref{thm:solvable} and Lemma~\ref{lem:fong}, there exist
		a finite solvable group $L$ and a normalized cocycle $\alpha$ such that
		\begin{equation}\label{eq:integral-model}
			\cO Gb\sim_{\Mor}\cO^\alpha L,
			\qquad O_{p'}(L)=1,
			\qquad \alpha\in Z^2(L,\mu_m(\cO)),\quad (m,p)=1.
		\end{equation}
		Let $\bar\alpha$ denote the reduction of $\alpha$ modulo the maximal ideal of $\cO$. Reducing the Morita equivalence gives
		\[
		k^{\bar\alpha}L\sim_{\Mor}kGb\sim_{\Mor}kP.
		\]
		Corollary~\ref{cor:twisted} implies that $L$ is a $p$-group.
		
		Since $(|L|,m)=1$, we have $H^2(L,\mu_m(\cO))=0$ and hence
		$\cO^\alpha L\cong\cO L$ \cite[Propositions~1.2.9(i) and~1.2.6]{Linckelmann2018a}.
		Thus \eqref{eq:integral-model} yields a Morita equivalence over $\cO$ between $\cO Gb$ and $\cO L$.
		Puig's theorem \cite[Theorem~8.2]{Puig1999} now implies that $b$ is nilpotent.
	\end{proof}
	
	\section*{Data availability}
	No datasets were generated or analysed in this work.
	
	\section*{Use of AI-assisted tools}
	During the preparation of this work, the authors used AI-assisted tools for literature searches, the exploration of possible proof strategies, the review of intermediate arguments, language editing, and LaTeX formatting. The authors carefully reviewed and edited all AI-generated content and take full responsibility for the content of this manuscript.

\end{document}